\documentclass[11pt,a4paper]{article}
\usepackage[T1]{fontenc}
\usepackage[utf8]{inputenc}
\usepackage[english]{babel}
\usepackage{amsmath,amssymb,amsthm,mathtools}
\usepackage[top=2.1cm,bottom=2.1cm,left=2.4cm,right=2.4cm]{geometry}
\usepackage{enumitem,booktabs,array,microtype,float}
\usepackage{xurl}
\usepackage[unicode,hidelinks]{hyperref}
\hypersetup{pdftitle={Extremal distributions for (a:b)-cuts in hypergraphs and the coupon collector's problem with group drawings},pdfkeywords={cuts, hypergraphs, coupon collector's problem, group drawings, optimal distributions},pdfauthor={N. I. Zderzhikov; D. A. Shabanov}}
\newtheorem{theorem}{Theorem}[section]
\newtheorem{lemma}[theorem]{Lemma}
\newtheorem{proposition}[theorem]{Proposition}

\theoremstyle{definition}

\theoremstyle{remark}
\newtheorem{remark}[theorem]{Remark}

\renewcommand{\P}{\mathbb P}
\newcommand{\E}{\mathbb E}
\newcommand{\1}{\mathbf 1}

\numberwithin{equation}{section}
\setlist[enumerate]{label=\textup{(\roman*)},leftmargin=*,itemsep=3pt}
\setlist[itemize]{leftmargin=*,itemsep=3pt}
\allowdisplaybreaks[2]
\title{Extremal distributions for $(a:b)$-cuts in hypergraphs and the coupon collector's problem with group drawings}
\author{N.\,I.~Zderzhikov\thanks{\raggedright National Research University Higher School of Economics, Moscow.
\mbox{E-mail: \href{mailto:nzderzhikov@hse.ru}{\nolinkurl{nzderzhikov@hse.ru}}}.}
\and D.\,A.~Shabanov\thanks{\raggedright Moscow Institute of Physics and Technology (National Research University) and National Research University Higher School of Economics, Moscow.
\mbox{E-mail: \href{mailto:shabanov.da@mipt.ru}{\nolinkurl{shabanov.da@mipt.ru}}}.}}
\date{}
\begin{document}
\maketitle
\begin{abstract}
In an independent $(a:b)$-coloring of a hypergraph, the color sets assigned
to the vertices are independent samples from a common distribution on
the $b$-element subsets of a palette of $a$ colors. An edge is cut if its vertices have
no color in common, or equivalently, if the complements of their color
sets cover the palette. In the coupon collector's problem with group
drawings, each such complement represents a draw of $s=a-b$ coupons.
Thus the maximum probability $\alpha_{\max}(a,b,k)$ of cutting an edge of
size $k$ equals the maximum probability of collecting all $a$ coupons
within $k$ independent draws.

When $b\mid a$, we give a short direct proof of the known formula for
$\alpha_{\max}$, determine the minimum expected collection time, and
characterize all optimal distributions for both criteria. For arbitrary
$a$ and $b$, we obtain restrictions on the sizes of the components of the
support of an optimal distribution. When $s\mid a$, choosing a block
uniformly at random from a fixed partition into blocks of size $s$
maximizes the probability of complete collection for every number of
draws and minimizes the expected collection time. For $b=2$ and for
$s=2$, we characterize all optimal distributions for both criteria.
When $s=2$ and $a$ is odd, for each distribution $y$ on drawn pairs
we can choose a distribution $y^*$ on the edges of a triangle and a
perfect matching on the remaining vertices such that collecting all
coupons within $k$ draws is at least as likely under $y^*$ as under $y$
for every $k$. The choice of $y^*$ is independent of $k$. In particular, we prove the
conjectures on the optimality of partitions when $s\mid a$ and on the
distributions that minimize the expected collection time when $s=2$.


\end{abstract}

\section{Introduction}\label{intro:section}

The maximum cut problem asks for a coloring of the vertices of a
hypergraph that cuts as many edges as possible. In the $(a:b)$-colorings
considered here, each vertex is assigned a $b$-element subset of a palette
of $a$ colors. An edge is cut if its vertices have no
color in common. For $b=1$ this is the usual maximum $a$-cut problem.

Suppose that the color sets assigned to the vertices are drawn
independently from a common distribution. Let $\alpha_{\max}(a,b,k)$
denote the largest probability, over all such distributions, that an
edge of size $k$ is cut.

The same parameter arises in the coupon collector's problem with group
drawings. If the color set $A\subset[a]$ assigned to a vertex is
interpreted as the set of missing coupons, then the set drawn is
$D=[a]\setminus A$, of size $s=a-b$, and
\[
 A_1\cap\cdots\cap A_k=\varnothing
 \quad\Longleftrightarrow\quad D_1\cup\cdots\cup D_k=[a].
\]
Let $x$ be the common distribution of the missing sets and $T_x$ the
time needed to collect all coupons. We write
\begin{equation}\label{intro:objectives}
 \Phi_k(x)\coloneqq\P(T_x\le k),\qquad
 \alpha_{\max}(a,b,k)\coloneqq\max_x\Phi_k(x),\qquad
 M(a,b)\coloneqq\min_x\E T_x.
\end{equation}
All optimizations are over distributions on sets of fixed size.
In terms of the drawn sets, $y_D=x_{[a]\setminus D}$, and we write
$T_y$ for the collection time and $\Phi_k(y)$ for the probability of
collecting all coupons within $k$ draws. Section~\ref{pre:model} gives the
precise formulation.

\subsection{Previous work}

Zderzhikov and Shabanov introduced $\alpha_{\max}$ in their study of the
maximum $(a:b)$-cut of a random hypergraph~\cite[Section~1.1]{ZderzhikovShabanov}.
The leading term in both the upper and the lower asymptotic bounds
on the cut size in~\cite[Theorems~1.1 and~1.2]{ZderzhikovShabanov}
is determined by $\alpha_{\max}$. The choice of a maximizing
distribution also matters for the refined lower bound. This bound
depends on the structure of the distribution and on the probabilities
that changing the color set of one vertex turns an uncut edge into a
cut edge or vice versa.

When $b=a-1$, we have $s=1$, so each draw yields a single coupon.
The uniform distribution on the $a$ coupons maximizes the probability
of complete collection for every number of draws. In particular, for $k\ge1$,
\begin{equation}\label{intro:classical}
 \begin{aligned}
 \alpha_{\max}(a,a-1,k)
 &=\sum_{j=0}^{a}(-1)^j\binom aj\left(1-\frac ja\right)^k,\\
 M(a,a-1)&=aH_a,\qquad H_a=\sum_{j=1}^a\frac1j.
 \end{aligned}
\end{equation}
Boneh and Hofri~\cite[Section~3, item~(9)]{BonehHofri1997} give a proof
of the optimality of the uniform distribution and attribute the argument
to S.~M.~Ross.

The formula for $\alpha_{\max}(a,a-1,k)$ also follows from Lemma~2
of Shabanov's paper on panchromatic colorings~\cite[pp.~621--622]{Shabanov2010Panchromatic}.
Suppose that $N$ vertices are colored with $a$ colors and $k$ vertices
are sampled uniformly without replacement. The probability that the
sample contains all $a$ colors is maximized when the sizes of the
color classes differ by at most one. For fixed $k$, as $N\to\infty$
with the class proportions tending to $p_1,\ldots,p_a$, the joint
distribution of the sampled colors converges to that of $k$ independent
draws with probabilities $p_1,\ldots,p_a$. For balanced color classes,
each limiting proportion is $1/a$. The probability of seeing all $a$
colors in $k$ independent draws is therefore maximized by the uniform
distribution, giving the stated value of $\alpha_{\max}(a,a-1,k)$
in~\eqref{intro:classical}.

Anceaume et al.~\cite{Anceaume2015,Anceaume2016} obtained more
general optimization results allowing for a prescribed number of
distinct coupons to be collected and for draws that yield no coupon.

Stadje~\cite{Stadje1990} studied the case in which sets of fixed size
are drawn uniformly at random. Adler and Ross~\cite{AdlerRoss2001}
considered the general model of independent random subsets with an
arbitrary common distribution. Falgas-Ravry, Larsson, and
Markström~\cite{FalgasRavry2020} gave general conditions for fast
coverage and concentration of the collection time. In these works
the distribution of the drawn set is given. When the set size is fixed,
one may also ask which distribution minimizes the expected collection
time.

Let $y^{\mathrm{all}}$ be the uniform distribution on all $s$-element
sets, so that $y^{\mathrm{all}}_D=\binom{a}{s}^{-1}$ for each
$D\in\binom{[a]}s$. When $s\mid a$, we can compare this distribution
with the uniform distribution on the blocks of a fixed partition into
$s$-element sets. The two distributions coincide for $s=1$ and differ
for $1<s<a$. Schilling~\cite{Schilling2021} studied the role of
$y^{\mathrm{all}}$ and conjectured that it is suboptimal for the
nontrivial range of package sizes. Berend and Sher~\cite{BerendSher2026}
proved the strict inequality
\begin{equation}\label{intro:nonuniform}
 \E T_{y^{\mathrm{all}}}>M(a,a-s),\qquad 2\le s\le a-2.
\end{equation}
This result rules out the uniform distribution on all packages, but
leaves open which distributions are optimal.

Two special cases of our problem follow directly from classical results.
For $k=2$, the probability of complete collection is a quadratic form
associated with the Kneser graph on the $b$-subsets of $[a]$, in which
two vertices are adjacent if their sets are disjoint. Its maximum and all cases of equality follow from the
Motzkin--Straus theorems~\cite{MotzkinStraus1965}. In particular,
\begin{equation}\label{intro:k2}
 \alpha_{\max}(a,b,2)=1-\frac1{\lfloor a/b\rfloor}.
\end{equation}
When $a<2b$, the right-hand side is zero, as it is impossible to
collect all coupons in two draws.

For $b\mid a$, the value of $\alpha_{\max}$ follows from Proposition~4
of Shabanov and Shaikheeva~\cite{ShabanovShaikheeva2020} on multiple
covers by independent sets. Passing from the sampling without
replacement considered there to independent draws from a common
distribution gives the universal bound
\begin{equation}\label{intro:bdiv-tail}
 \P(A_1\cap\cdots\cap A_k\ne\varnothing)
 \ge\left(\frac ba\right)^{k-1},\qquad k\ge1.
\end{equation}
When $a=qb$, the bound is attained by drawing uniformly from the
$b$-element blocks of a fixed partition. Hence
\begin{equation}\label{intro:bdiv-probability}
 \alpha_{\max}(qb,b,k)=1-q^{1-k},\qquad q\ge2,\quad k\ge1.
\end{equation}

Fang and Chang~\cite{FangChang2014} solved the problem for packages
containing half of all coupons. When $a=2s$, the uniform distribution
on two complementary sets is optimal for each number of draws and
therefore also for the expected collection time:
\begin{equation}\label{intro:half}
 \alpha_{\max}(2s,s,k)=1-2^{1-k}\quad(k\ge1),
 \qquad M(2s,s)=3.
\end{equation}
For any divisor $s$ of $a$, a natural candidate is the uniform
distribution on the blocks of a partition into $a/s$ sets of size $s$.
Fang and Chang~\cite{FangChang2014} discussed the conjecture that
such a partition is optimal. Schilling's
dissertation~\cite[p.~72, Conjecture~6.2.1]{Schilling2020} states the conjecture for all $s\ge3$ with $s\mid a$ and asserts that
the uniform distributions on the blocks of such partitions are the only
minimizers of the expected collection time. Theorems~\ref{sdiv:thm:main} and~\ref{sdiv:cor:mean}
prove this conjecture for sets of fixed size $s$. Moreover, drawing
uniformly from the blocks of a partition maximizes the probability of
complete collection for every number of draws.

Caron, Hlynka, and McDonald~\cite{CHM1992} studied the case $s=2$
in their analysis of algorithms for detecting necessary linear
constraints. They considered distributions that minimize the expected
collection time. For an even number of coupons they proposed the uniform
distribution on a perfect matching. For an odd number they proposed a
distribution supported on the edges of a triangle and a perfect matching
on the remaining vertices. These supports may be written as
\begin{equation}\label{intro:pair-supports}
 \begin{aligned}
 a=2q:\quad&qK_2, &&y_e=1/q,\\
 a=2m+3:\quad&K_3\sqcup mK_2,
 &&y_e=\begin{cases}
 v,&e\in E(K_3),\\
 u,&e\in E(mK_2),
 \end{cases}
 \end{aligned}
\end{equation}
where $K_2$ is a single edge, $K_3$ is a triangle, and $\sqcup$
denotes disjoint union. In the odd case the weights satisfy
$3v+mu=1$. The authors of~\cite{CHM1992} established global optimality
in the small cases $a=2,3,4$ and obtained results on strict local
minima for the corresponding families at larger values of $a$. They
conjectured global optimality for an arbitrary number of coupons.
Schilling~\cite[p.~28, Conjectures~3.2.1 and~3.2.2]{Schilling2020}
formulated the conjecture with a characterization of all minimizers:
for even $a$, the uniform distribution on a perfect matching, and for
odd $a\ge5$, a triangle and isolated edges with optimal weights.
Theorem~\ref{s2:thm:mean} confirms this conjecture and characterizes
all minimizers of the expected collection time.

\subsection{Main results and outline}

In each of the four cases $b\mid a$, $s\mid a$, $b=2$, and $s=2$,
we determine both the maximum probability $\alpha_{\max}(a,b,k)$ of collecting all coupons
within $k$ draws and the minimum expected collection time $M(a,b)$.
We characterize all optimal distributions for each objective. For
$k=2$ and arbitrary $1\le b<a$, we also determine $\alpha_{\max}(a,b,2)$
and all distributions attaining it.

The \emph{support of a distribution} is the family of sets of positive
probability. These are missing sets for $x$ and drawn sets for $y$.
We regard them as the edges of a hypergraph on the set of coupons
$[a]$. We call its connected components the \emph{components of the
support}. The size of a component is its number of vertices. Each
coupon belonging to no edge forms a component of size~$1$.
Table~\ref{intro:table} lists the supports up to a permutation of the
coupons.

\begin{table}[H]
\centering
\small
\renewcommand{\arraystretch}{1.2}
\begin{tabular}{@{}l@{\hspace{1.5em}}>{\raggedright\arraybackslash}p{0.73\linewidth}@{\hspace{1.5em}}c@{}}
\toprule
Case & Support of an optimal distribution & Section\\
\midrule
$k=2$ & A complete $\lfloor a/b\rfloor$-partite induced subgraph of
the Kneser graph on the $b$-subsets of $[a]$ & \ref{k2:section}\\
$b\mid a$ & A partition into missing $b$-sets & \ref{bdiv:section}\\
$s\mid a$ & A partition into drawn $s$-sets & \ref{sdiv:section}\\
$b=2$ & For even $a$, a perfect matching. For odd $a$, a triangle and
a perfect matching on the remaining vertices & \ref{b2:section}\\
$s=2$ & For even $a$, a perfect matching. For odd $a$, a triangle and
a perfect matching on the remaining vertices & \ref{s2:sec:main}\\
\bottomrule
\end{tabular}
\caption{Supports of optimal distributions. In the last two rows the
edges represent missing and drawn two-element sets, respectively.
The degenerate cases and additional cases of equality are given in
the corresponding sections.}
\label{intro:table}
\end{table}

\paragraph{Missing sets when $b\mid a$ and components of the support.}
Section~\ref{bdiv:section} gives a short direct proof
of~\eqref{intro:bdiv-probability} for $a=qb$ by bounding the probability
that the missing sets have a nonempty intersection in terms of the
expected size of that intersection. The equality conditions show that,
for each $k\ge2$, a distribution maximizes $\Phi_k$ if and only if
the missing set is chosen uniformly from the $q$ blocks of a
partition (Theorem~\ref{bdiv:main}). Such a distribution is optimal
simultaneously for all $k$. Summing the tail probabilities gives
\[
 M(qb,b)=2+\frac{1}{q-1}.
\]
These distributions are also the only minimizers of the expected
collection time (Theorem~\ref{bdiv:mean}).

For arbitrary $1\le b<a$, the partition result applies to a union of
components of the support whenever its number of vertices is divisible
by $b$. Suppose that a distribution maximizes $\Phi_k$ for some $k\ge2$
or minimizes the expected collection time. Then in any such nonempty
union of components of the hypergraph of missing sets, each component
has exactly $b$ vertices. In particular, the number of components whose
sizes are not divisible by $b$ is at most $b-\gcd(a,b)$
(Theorem~\ref{bdiv:components}).

\paragraph{Optimality of uniform sampling from a partition when $s\mid a$.}
Let $a=qs$, where $q\ge2$, and fix a partition into $q$ blocks of size
$s$. Under the distribution $y^{\mathrm{part}}$, each draw selects one
of these blocks uniformly at random. Theorem~\ref{sdiv:thm:main} gives
\begin{equation}\label{intro:sdiv-domination}
 \forall y\ \forall k\ge0:\qquad
 \Phi_k(y)\le\Phi_k(y^{\mathrm{part}}).
\end{equation}
This yields the explicit values
\begin{equation}\label{intro:sdiv-results}
 \begin{aligned}
 \alpha_{\max}(qs,(q-1)s,k)
 &=\sum_{j=0}^{q}(-1)^j\binom qj\left(1-\frac jq\right)^k,\\
 M(qs,(q-1)s)&=qH_q.
 \end{aligned}
\end{equation}
For $k=0$ we use the convention $0^0=1$. For $k\ge q$ the uniform
distributions on the blocks of partitions are the only maximizers of
the probability of complete collection. They are also the only
minimizers of the expected collection time. For $k<q$, complete
collection is impossible under any distribution. Thus the uniform
distribution on the blocks of any fixed partition into $s$-sets is
optimal simultaneously for all $k$.

For arbitrary $1\le s<a$, Theorem~\ref{sdiv:cor:components} gives an
analogue of Theorem~\ref{bdiv:components} for the hypergraph of drawn
sets. Suppose that a distribution maximizes $\Phi_k$ for some
$k\ge\lceil a/s\rceil$ or minimizes the expected collection time.
If the total number of vertices in a nonempty union of components
is divisible by $s$, then each of these components has exactly $s$
vertices. The number of components whose sizes are not divisible by
$s$ is at most $s-\gcd(a,s)$.

\paragraph{Collection probability and expected collection time for $b=2$.}
For even $a$ the canonical support is a perfect matching. For odd $a$
it consists of a triangle and isolated edges. Each edge represents a
missing set. For odd $a=2m+3$, $m\ge1$, we obtain the explicit formula
\begin{equation}\label{intro:b2-probability}
 \alpha_{\max}(2m+3,2,k)
 =1-\frac{2^k-1}{(m\rho_k+3)^{k-1}},
 \qquad \rho_k=(2^k-1)^{1/(k-1)},\quad k\ge2.
\end{equation}
For $k\ge3$, every maximizer has support $K_3\sqcup mK_2$, with
uniquely determined optimal weights. For $k=2$ there may be additional
cases of equality. Minimizing the expected collection time also yields
a unique choice of weights on the canonical support
(Theorem~\ref{b2:cor:mean}).

\paragraph{A dominating family for $s=2$.}
Let $\mathcal C_a$ denote the family of canonical distributions on
drawn pairs. For even $a$ these are the uniform distributions on the
edges of perfect matchings. For $a=2m+3$ they are distributions on
$K_3\sqcup mK_2$ with a common weight for the triangle edges and a
common weight for the isolated edges. Theorem~\ref{s2:thm:dom} states
that
\begin{equation}\label{intro:s2-domination}
 \forall y\ \exists y^*\in\mathcal C_a\ \forall k\ge0:
 \qquad \Phi_k(y)\le\Phi_k(y^*).
\end{equation}
For each $y$, a single distribution $y^*$ with fixed weights works
for every $k$. The weights of a maximizer of $\Phi_k$ may depend
on $k$. The tail-sum formula~\eqref{pre:tail-equation} gives
$\E T_{y^*}\le\E T_y$, so both optimization problems reduce to the
family $\mathcal C_a$:
\begin{equation}\label{intro:s2-results}
 \alpha_{\max}(a,a-2,k)=\max_{y\in\mathcal C_a}\Phi_k(y),
 \qquad M(a,a-2)=\min_{y\in\mathcal C_a}\E T_y.
\end{equation}
When $a=2m+3\ge5$, this family has a single parameter, the total
probability mass of the triangle. Complete collection is possible
starting at $k=m+2$. For each $k\ge m+2$ the mass maximizing the collection probability is unique.
The mass minimizing the expected collection time is also unique
(Theorems~\ref{s2:thm:horizon} and~\ref{s2:thm:mean}). When the number of coupons is odd and the number of draws is the
smallest for which complete collection is possible, we obtain
\begin{equation}\label{intro:s2-first}
 \alpha_{\max}(2k-1,2k-3,k)=\frac{4k!}{3k^k},\qquad k\ge2.
\end{equation}

\paragraph{Main ideas of the proofs.}
When $a=qs$, we construct an injection from covering sequences of
numbered tickets bearing coupon sets to sequences of colored tickets in
which all $q$ colors occur. Several numbered tickets may bear the same
set. For any distribution on drawn sets, the construction gives a single
distribution on the $q$ colors for which the comparison
in~\eqref{sdiv:eq:comparison} holds for every $k$. It remains to
equalize the color probabilities and analyze the cases of equality.

For $b=2$, write $x_{ij}$ for the probability that $\{i,j\}$ is the
missing pair. The probability of incomplete collection and the expected
collection time are given in terms of the same difference of sums,
\[
 \sum_i\varphi(p_i)-\sum_{i<j}\varphi(x_{ij}),
 \qquad p_i=\sum_{j\ne i}x_{ij}.
\]
For $\varphi(t)=t^k$, this difference equals $1-\Phi_k(x)$ for $k\ge1$.
For $\varphi(t)=t^2/(1-t)$, it equals $\E T_x-2$ whenever the
expectation is finite. Transfers of mass and variations that preserve
the sum of the weights at each vertex reduce both problems to
canonical supports.

For $s=2$, we isolate edges of the support graph, vary the weights on $4$-cycles, and
transform a bouquet of triangles. These operations rule out
noncanonical supports both for distributions maximizing the collection
probability when $k\ge\lceil a/2\rceil$ and for those minimizing the
expected collection time.

On the support $K_3\sqcup mK_2$, we express the expected collection
time in terms of the time until all isolated edges have appeared and
the number of additional draws needed to cover the triangle. By
Lemma~\ref{s2:lem:concavity}, the collection probabilities are concave
in the mass of the triangle on an interval common to all $k\ge m+2$.
This concavity and the signs of the derivatives imply that the
expected collection time is minimized at a unique value of the mass.
We also use this concavity and the characterization of the maximizers
of finite nonnegative linear combinations of the collection
probabilities (Theorem~\ref{s2:thm:structure}) to obtain a simultaneous
comparison for all horizons.

\section{Two formulations and a common extremal parameter}\label{pre:model}

\subsection{Cuts under independent assignments of types}
Let $a,b$ be integers with $1\le b<a$, and let $s=a-b$.
Write $[a]=\{1,\ldots,a\}$ and $\mathcal A=\binom{[a]}b$.
We call the elements of $\mathcal A$ types and denote by
$\Delta_{a,b}$ the simplex of distributions
\[
 x=(x_A)_{A\in\mathcal A},\qquad x_A\ge0,\qquad\sum_Ax_A=1.
\]
For independent types $A_1,\ldots,A_k$ with common distribution $x$,
define
\begin{equation}\label{pre:alpha}
 \Phi_k(x)\coloneqq\P(A_1\cap\cdots\cap A_k=\varnothing),\qquad
 \alpha_{\max}(a,b,k)\coloneqq\max_{x\in\Delta_{a,b}}\Phi_k(x).
\end{equation}
For $k=0$, we take the intersection of the empty family to be $[a]$.
In particular, $\Phi_0(x)=\Phi_1(x)=0$.

For a $k$-uniform hypergraph $H=(V,E)$ and a map
$\sigma:V\to\mathcal A$, let
\[
 C_H(\sigma)=\bigl|\{e\in E:\textstyle\bigcap_{v\in e}\sigma(v)=\varnothing\}\bigr|,
 \qquad \mu_{a:b}(H)=\max_\sigma C_H(\sigma).
\]
All hypergraphs are finite and have no multiple edges. When types
are assigned independently according to $x$, each edge of size $k$
is cut with probability $\Phi_k(x)$. By linearity of expectation,
\begin{equation}\label{pre:cut-bound}
 \E C_H(\sigma)=\Phi_k(x)|E(H)|,
 \qquad
 \mu_{a:b}(H)\ge\alpha_{\max}(a,b,k)|E(H)|.
\end{equation}
The maximum in~\eqref{pre:alpha} is over the common distribution of
the independently chosen vertex types. In the definition of
$\mu_{a:b}(H)$ the maximum is over all colorings of the fixed
hypergraph $H$. Thus $\alpha_{\max}(a,b,k)$ gives a guaranteed
fraction of edges that can be cut. For example, when $a=qs$,
$b=(q-1)s$, and $q\ge2$, Theorem~\ref{sdiv:thm:main} gives, for $k\ge1$,
\[
 \mu_{qs:(q-1)s}(H)\ge
 \left[\sum_{j=0}^{q}(-1)^j\binom qj
       \left(1-\frac jq\right)^k\right]|E(H)|.
\]
Sharper asymptotic bounds on the cut size in random hypergraphs,
including a refined lower bound that takes the optimal type
distribution into account, are given
in~\cite[Theorems~1.1 and~1.2]{ZderzhikovShabanov}.

\subsection{The time to collect all coupons}
At each step we independently choose a missing set $A_t$ from the
distribution $x\in\Delta_{a,b}$, and the collector receives its
complement $D_t=[a]\setminus A_t$, of size $s$. Define the
collection time by
\[
 T_x\coloneqq\inf\{k\ge0:D_1\cup\cdots\cup D_k=[a]\},
 \qquad \inf\varnothing=\infty.
\]
By De Morgan's law, $\Phi_k(x)=\P(T_x\le k)$. Write $y$ for the
distribution of the drawn set, so that $y_D=x_{[a]\setminus D}$.
In this notation we write $T_y$ and $\Phi_k(y)$. Define the
minimum expected collection time by
\begin{equation}\label{pre:mean}
 M(a,b)\coloneqq\min_{x\in\Delta_{a,b}}\E T_x.
\end{equation}
The parameters $a$ and $b$ denote the total number of coupons and the
size of the \emph{missing} set, respectively. All statements about
uniqueness of distributions are understood up to a permutation of
the coupons.

\begin{proposition}[Attainment of the extrema and the tail-sum formula]\label{pre:tail}
The maximum in~\eqref{pre:alpha} and the minimum in~\eqref{pre:mean}
are attained. For every distribution $x$, the following identity
holds in $[0,\infty]$:
\begin{equation}\label{pre:tail-equation}
 \E T_x=\sum_{k=0}^{\infty}\P(T_x>k)
       =\sum_{k=0}^{\infty}(1-\Phi_k(x)).
\end{equation}
The expected collection time is finite if and only if each coupon
has positive probability of being drawn.
\end{proposition}
\begin{proof}
The function $\Phi_k$ is a polynomial in the coordinates of $x$ and
therefore attains its maximum on the compact simplex.
Formula~\eqref{pre:tail-equation} follows from
$T_x=\sum_{k\ge0}\1_{\{T_x>k\}}$ and the monotone convergence theorem.
Write $p_i=\P(i\in A_t)$. If $p_i=1$ for some $i$, complete
collection is impossible. Otherwise,
$\P(T_x>k)\le\sum_{i=1}^a p_i^k$, which gives $\E T_x<\infty$.
The partial sums in~\eqref{pre:tail-equation} are continuous and
nondecreasing in the number of terms, so $x\mapsto\E T_x$ is lower semicontinuous
and attains its minimum on the simplex. The expected collection time
is finite under the uniform distribution on all drawn $s$-sets.
\end{proof}

\begin{proposition}[Comparison of expected collection times]\label{pre:transfer}
Suppose that a distribution $x_*$ satisfies
$\Phi_k(x)\le\Phi_k(x_*)$ for all $x\in\Delta_{a,b}$ and all $k\ge0$.
Then $\E T_{x_*}=M(a,b)$. If for every distribution $x$ outside a
family $\mathcal F\subseteq\Delta_{a,b}$ there is some $k$ for which
$\Phi_k(x)<\Phi_k(x_*)$, then every minimizer of the expected collection
time belongs to $\mathcal F$.
More generally, suppose that for every $x\in\Delta_{a,b}$ there is a
distribution $z=z(x)\in\mathcal C$ such that $\Phi_k(x)\le\Phi_k(z)$
simultaneously for all $k\ge0$. Then $\mathcal C$ contains at least one
minimizer of the expected collection time.
\end{proposition}
\begin{proof}
Summing the inequalities for the tail probabilities
in~\eqref{pre:tail-equation} gives $\E T_{x_*}\le\E T_x$ for every $x$.
The minimum is finite by Proposition~\ref{pre:tail}. Hence, if the
inequality is strict for at least one tail probability, then
$\E T_{x_*}<\E T_x$. In the last part of the proposition, take $x$ to
be a minimizer of the expected collection time. Then
$\E T_{z(x)}\le\E T_x=M(a,b)$, so $z(x)$ also attains the minimum.
\end{proof}

\section{Collection in two draws: the case \texorpdfstring{$k=2$}{k=2}}\label{k2:section}

When $k=2$, the problem reduces to maximizing the quadratic form of the
Kneser graph on the $b$-subsets of $[a]$. The results of
Motzkin and Straus~\cite{MotzkinStraus1965} applied to this graph give
the value of $\alpha_{\max}(a,b,2)$ and all cases of equality.

\begin{theorem}\label{k2:main}
Let $1\le b<a$ and $q=\lfloor a/b\rfloor$. Then
\begin{equation}\label{k2:value}
 \alpha_{\max}(a,b,2)=1-\frac1q.
\end{equation}
A distribution $x$ attains the maximum if and only if its support admits
a partition $\mathcal S=\mathcal A_1\sqcup\cdots\sqcup\mathcal A_q$
into nonempty families satisfying
\begin{equation}\label{k2:equality}
\begin{aligned}
 A\cap B&\ne\varnothing
 &&(A,B\in\mathcal A_i),\\
 A\cap B&=\varnothing
 &&(A\in\mathcal A_i,\ B\in\mathcal A_j,\ i\ne j),\\
 \sum_{A\in\mathcal A_i}x_A&=\frac1q
 &&(i\in[q]).
\end{aligned}
\end{equation}
\end{theorem}

\begin{proof}
Two drawn sets cover $[a]$ if and only if the corresponding missing sets
are disjoint. Hence
\begin{equation}\label{k2:quadratic}
 \Phi_2(x)=
 \sum_{\substack{A,B\in\binom{[a]}b\\A\cap B=\varnothing}}x_Ax_B,
\end{equation}
where the sum is over ordered pairs.

By Proposition~\ref{pre:tail}, the maximum is attained. Choose a maximizer
$x$ with support $\mathcal S$. We show that the relation $A\cap B\ne\varnothing$
is transitive on $\mathcal S$. Suppose to the contrary that $A\cap B\ne\varnothing$,
$B\cap C\ne\varnothing$, and $A\cap C=\varnothing$. For sufficiently
small $\varepsilon>0$, both perturbations
\[
 x_A^{\pm}=x_A\pm\varepsilon,\qquad
 x_C^{\pm}=x_C\pm\varepsilon,\qquad
 x_B^{\pm}=x_B\mp2\varepsilon
\]
preserve the positivity of the weights and their sum. Of the three sets
whose weights change, only $A$ and $C$ are disjoint. The linear terms
cancel when the two increments are added, giving
\[
 \Phi_2(x^+)+\Phi_2(x^-)-2\Phi_2(x)=4\varepsilon^2>0.
\]
At least one perturbation increases $\Phi_2$, a contradiction.

This relation is also reflexive and symmetric, so it is an equivalence relation on $\mathcal S$.
Let $\mathcal A_1,\ldots,\mathcal A_r$ be its equivalence classes and
let $w_i=\sum_{A\in\mathcal A_i}x_A$. Choosing one set from each class
gives $r$ pairwise disjoint $b$-subsets, so $r\le q$. The collection is
complete precisely when the two selected missing sets belong to different
classes. Since $\sum_i w_i=1$, we obtain
\[
 \Phi_2(x)=1-\sum_{i=1}^r w_i^2
 \le 1-\frac1r\le 1-\frac1q.
\]
The uniform distribution on $q$ pairwise disjoint $b$-subsets gives
$1-1/q$. This is therefore the maximum, and both inequalities are
equalities for the chosen $x$. Hence $r=q$ and $w_i=1/q$ for each $i$.
Conversely, every distribution satisfying~\eqref{k2:equality} gives the
same probability of complete collection.
\end{proof}

\begin{remark}\label{k2:degenerate}
If $a<2b$, then $q=1$. Any two missing sets intersect, so
$\alpha_{\max}(a,b,2)=0$ and every distribution is optimal.
\end{remark}

\section{Missing sets when \texorpdfstring{$b\mid a$}{b | a}}\label{bdiv:section}

When $b\mid a$, the value of $\alpha_{\max}$ follows from
Proposition~4 of Shabanov and Shaikheeva~\cite{ShabanovShaikheeva2020}
by passing from sampling without replacement to independent draws.
The equality conditions in the bound on the probability of a nonempty
intersection determine all maximizers of $\Phi_k$ for $k\ge2$.
Since each is optimal simultaneously for all $k$, the same bound also
determines the minimizers of the expected collection time.

\begin{theorem}[Probability of complete collection and optimal distributions]\label{bdiv:main}
Let $a=qb$, where $q\ge2$ is an integer. Then for every $k\ge1$,
\begin{equation}\label{bdiv:value}
 \alpha_{\max}(qb,b,k)=1-q^{1-k}.
\end{equation}
For each fixed $k\ge2$, the maximum is attained precisely by the uniform
distributions on the blocks of a partition of $[a]$ into $q$ sets of
size $b$. Each such distribution is optimal simultaneously for all $k$.
\end{theorem}

The main estimate compares the probability that the intersection is
nonempty with its expected size. A nonempty intersection contains at most
$b$ elements, giving the first inequality. The second follows by averaging
the probabilities that individual coupons are missing.

\begin{lemma}[A bound on the probability of a nonempty intersection]\label{bdiv:tail}
For any distribution $x$, any $1\le b<a$, and any $k\ge1$, let
$p_i=\P(i\in A_1)$. Then
\begin{equation}\label{bdiv:bound}
 \P(T_x>k)\ge\frac1b\sum_{i=1}^a p_i^k
 \ge\left(\frac ba\right)^{k-1}.
\end{equation}
\end{lemma}

\begin{proof}
Let $Z_k=|A_1\cap\cdots\cap A_k|$. Then
$Z_k\le b\,\mathbf1_{\{T_x>k\}}$. By independence of the draws,
$\E Z_k=\sum_i p_i^k$, which gives the first inequality.
Since $\sum_i p_i=b$, Jensen's inequality gives
\[
 \frac1b\sum_{i=1}^a p_i^k
 \ge\frac ab\left(\frac1a\sum_{i=1}^a p_i\right)^k
 =\left(\frac ba\right)^{k-1}.\qedhere
\]
\end{proof}

\begin{proof}[Proof of Theorem~\ref{bdiv:main}]
Lemma~\ref{bdiv:tail} gives the required upper bound on $\Phi_k(x)$.
If the missing set is chosen uniformly from the blocks of a partition,
then for $k\ge1$ the intersection is nonempty precisely when the same
block is chosen on all $k$ draws. This event has probability
$q(1/q)^k=q^{1-k}$. Thus the bound is attained.

Now suppose equality holds for some $k\ge2$.
As in the proof of the lemma, let $Z_k=|A_1\cap\cdots\cap A_k|$.
Both inequalities in~\eqref{bdiv:bound} must be equalities. Equality in
the first requires $Z_k$ to be either $0$ or $b$ almost surely. In the
second, strict convexity of $u\mapsto u^k$ requires all the $p_i$ to
be equal. Thus
\[
 \P(0<Z_k<b)=0,\qquad p_i=\frac{b}{a}=\frac{1}{q}\quad(i\in[a]).
\]
If two distinct sets $G,H$ of positive weight intersect, then the event
$A_1=\cdots=A_{k-1}=G$, $A_k=H$ has positive probability. On this
event, $0<Z_k=|G\cap H|<b$, a contradiction.

Thus the sets of positive weight are pairwise disjoint. Since $p_i>0$
for every $i$, they cover $[a]$ and form a partition into $q$ blocks.
For $i\in G$ we have $x_G=p_i=1/q$.
\end{proof}

\subsection{Expected collection time}

\begin{theorem}[Minimum expected collection time]\label{bdiv:mean}
If $a=qb$ and $q\ge2$, then
\begin{equation}\label{bdiv:meanvalue}
 M(qb,b)=2+\frac1{q-1}=\frac{2a-b}{a-b}.
\end{equation}
The expected collection time is minimized precisely when the missing set is
chosen uniformly from the blocks of a partition of $[a]$ into $q$ sets
of size $b$.
\end{theorem}

\begin{proof}
Let $T_*$ denote the collection time for one of these partitions.
By Theorem~\ref{bdiv:main}, for every $k\ge0$,
$\P(T_x>k)\ge\P(T_*>k)$. The tail-sum formula therefore gives
\[
 \E T_x\ge\E T_*
 =1+\sum_{k=1}^{\infty}q^{1-k}
 =2+\frac1{q-1}.
\]
If $x$ is not uniform on the blocks of any partition of $[a]$ into
$b$-sets, the tail bound is strict at $k=2$. Hence the inequality for
the expected collection time is strict as well.
\end{proof}

\subsection{Restrictions on the components of an optimal support}

Let $a$ and $b$ now be arbitrary integers with $1\le b<a$.
Consider the hypergraph of missing sets of positive weight. If a union
of its components has positive mass and contains $mb$ vertices for some
$m\ge2$, Theorem~\ref{bdiv:main} applies to the conditional distribution
on this union.

\begin{theorem}[Components of the support of the missing sets]\label{bdiv:components}
Suppose $x$ attains $\alpha_{\max}(a,b,k)$ for some $k\ge2$ or
minimizes $\E T_x$. For any distinct components $C_1,\ldots,C_\ell$
of the support, where $\ell\ge1$,
\begin{equation}\label{bdiv:componentscondition}
 b\mid\bigl(|C_1|+\cdots+|C_\ell|\bigr)
 \quad\Longleftrightarrow\quad
 |C_1|=\cdots=|C_\ell|=b.
\end{equation}
There are at most $b-\gcd(a,b)$ components whose sizes are not divisible
by $b$. The sum of the sizes of any nonempty subfamily of these
components is also not divisible by $b$.
\end{theorem}

\begin{proof}
The intersection $A_1\cap\cdots\cap A_j$ can be nonempty only if all
the selected sets belong to the same component. Thus for $j\ge1$,
\begin{equation}\label{bdiv:decomposition}
 \P(T_x>j)=\sum_C
 \P(C\cap A_1\cap\cdots\cap A_j\ne\varnothing).
\end{equation}

Let $U$ be a union of components with $|U|=mb$, and suppose that $U$ contains an edge.
If $m=1$, this edge is $U$. Suppose $m\ge2$ and write
$\sigma=\sum_{A\subseteq U}x_A>0$. Conditional on all $j$ sets being
chosen from $U$, they are independent with a common distribution on the
$b$-subsets of $U$. By Theorem~\ref{bdiv:main}, for $j\ge2$,
\[
 \P(U\cap A_1\cap\cdots\cap A_j\ne\varnothing)
 \ge\sigma^j m^{1-j},
\]
with equality only if the conditional distribution is uniform on the
blocks of a partition of $U$.

Distribute the mass $\sigma$ equally among $m$ disjoint $b$-element
blocks covering $U$, keeping all weights outside $U$ fixed. The terms in
\eqref{bdiv:decomposition} for components outside $U$ stay the same.
If the original conditional distribution was not uniform on the blocks
of a partition, then the probability of incomplete collection strictly
decreases for every $j\ge2$. The new expected collection time is finite,
since it suffices to wait for two fixed distinct blocks, each chosen
with probability $\sigma/m>0$. By the tail-sum formula, the mean
collection time also strictly decreases. This contradicts optimality
under either criterion. Hence the original conditional distribution is
uniform on the blocks of a partition of $U$, and every component in $U$
has size $b$.

We next show that there are fewer than $b$ isolated vertices. Otherwise
we could add between one and $b$ isolated vertices to any component
containing an edge to obtain a union of size divisible by $b$. By the
preceding argument its edges would partition all of its vertices, which
is impossible in the presence of an isolated vertex. Thus every nonempty
union of components whose total size is divisible by $b$ contains an
edge. The preceding argument gives the forward implication
in~\eqref{bdiv:componentscondition}. The converse is immediate.

Let $g=\gcd(a,b)$, and let $t$ be the number of components whose sizes
are not divisible by $b$. The bound is immediate if $t=0$. The sum of the
sizes of these $t$ components is congruent to $a$ modulo $b$, so it is
divisible by $g$. Partition them into groups whose total sizes are each
divisible by $g$. For the first group choose an inclusion-minimal
nonempty subfamily with this property. Such a subfamily exists, since
the family of all $t$ components has the required property. After the
chosen group is removed, the remaining total size is still divisible
by $g$. Repeat until no components remain.

Suppose one of the chosen groups contains $h$ components. Order them
arbitrarily and consider the sums of the first $0,1,\ldots,h-1$ sizes.
If two sums are congruent modulo $g$, the components between the two
positions form a nonempty proper subfamily whose total size is divisible
by $g$. This contradicts minimality. The residues are therefore distinct,
and $h\le g$.

Suppose there are $m$ groups. The total sizes of the first
$0,1,\ldots,m$ groups are pairwise distinct modulo $b$. Otherwise some
nonempty union of groups would have total size divisible by $b$, contrary
to~\eqref{bdiv:componentscondition}. Each partial sum is divisible by
$g$, so its residue modulo $b$ belongs to $\{0,g,\ldots,b-g\}$.
There are $b/g$ such residues, and each group contains at most $g$
components. Hence
\[
 m+1\le\frac{b}{g},\qquad t\le gm\le b-g.\qedhere
\]
\end{proof}

The bound in Theorem~\ref{bdiv:components} can be attained even in
the absence of isolated vertices. Let $a=8$ and $b=3$, and split the
coupons into two sets of size four. On each draw, choose the missing triple uniformly from eight possibilities,
four in each set of four. Two selected triples are disjoint precisely when they come from
different sets of four. Hence $\Phi_2=1/2$, which equals
$\alpha_{\max}(8,3,2)$ by Theorem~\ref{k2:main}. The support consists
of two components of size four, and their number is $3-\gcd(8,3)=2$.

\begin{samepage}
\section{Drawn sets when \texorpdfstring{$s\mid a$}{s | a}}\label{sdiv:section}

Write $y_G=x_{[a]\setminus G}$ for the probability of drawing
$G\in\binom{[a]}s$, and let $T_y$ denote the collection time.

\begin{theorem}[Probability of complete collection and maximizers when $s\mid a$]\label{sdiv:thm:main}
Let $a=qs$, where $q\ge2$ is an integer. For every integer $k\ge0$,
\begin{equation}\label{sdiv:eq:maximum}
 \alpha_{\max}(qs,(q-1)s,k)
 =\sum_{j=0}^q(-1)^j\binom qj\left(1-\frac{j}{q}\right)^k ,
\end{equation}
where $0^0=1$. For $k\ge q$, the maximum is attained precisely when
$D_t$ is chosen uniformly from the blocks of a partition of $[a]$ into
$q$ sets of size $s$. For $k<q$, the probability of complete collection
is zero for every distribution.
\end{theorem}
\end{samepage}

\paragraph{Attainment of the bound.}
Partition $[a]$ into $q$ blocks of size $s$ and choose each block with
probability $1/q$. To collect all coupons, each block must appear at
least once. The probability of missing $j$ specified blocks is
$(1-j/q)^k$. Inclusion--exclusion gives the right-hand side of
\eqref{sdiv:eq:maximum}. For $k<q$, any $k$ sets cover at most $ks<a$
elements. In what follows, assume $k\ge q$.

\subsection{Comparison with the uniform distribution on the blocks of a partition}

Consider two bags containing the same number of tickets. Each ticket
in the first bag bears the set of coupons to be drawn. Each ticket in
the second bag has a color and a number within that color. We choose
the number of tickets of each of the $q=a/s$ colors using the
distribution $y$ of the drawn sets. To each sequence of $k$ tickets from
the first bag whose sets cover $[a]$, we assign a sequence of $k$ tickets
from the second bag in which all $q$ colors appear. Distinct covering
sequences will give distinct sequences of colored tickets. When tickets
are drawn uniformly with replacement, all sequences of a given length
are equally likely. Thus the probability of seeing all colors in $k$
draws will be at least $\Phi_k(y)$. We first prove the comparison for
rational probabilities, then extend it to arbitrary probabilities.

\subsubsection{Tickets bearing sets and colored tickets}

Let $p_i=\P(i\in D_t)$. Then
\begin{equation}\label{sdiv:eq:marginal}
 p_i=\sum_{G\ni i}y_G,\qquad \sum_{i=1}^a p_i=\sum_G |G|y_G=s.
\end{equation}
Number the elements so that $p_1\le\cdots\le p_a$.

First suppose that all $y_G$ are rational and let $N$ be a common
denominator. The first bag contains $N$ tickets numbered from $1$ to
$N$. Each ticket bears a set $G$, called its \emph{type}. There are
$m_G=Ny_G$ tickets of type $G$. A ticket drawn uniformly at random
has type $G$ with probability $m_G/N=y_G$. A sequence of $k$ tickets
drawn with replacement is called \emph{covering} if the union of their
types is $[a]$.

The second bag also contains $N$ tickets. Ticket $(j,r)$ has color $j$
and number $r$ within that color. Give color $j$ to $M_j$ tickets,
numbered from $1$ to $M_j$, where
\begin{equation}\label{sdiv:eq:colors}
 M_j=Np_{(j-1)s+1}\quad(1\le j<q),\qquad
 M_q=N-\sum_{j=1}^{q-1}M_j.
\end{equation}
The counts $M_j$ are fixed for all input sequences. We construct the
matching between the tickets in the two bags separately for each
covering sequence.

\subsubsection{Construction of the map}

Fix an order on the types. Let $z=(z_1,\ldots,z_k)$ be a covering
sequence of tickets from the first bag. At each of $q$ steps, we choose
one type $G_j$ that occurs in $z$ and assign tickets of color $j$.

\smallskip
\noindent\textit{Step $j$.}
Take the least element not yet covered by the selected types,
\begin{equation}\label{sdiv:eq:pivot}
 i_j=\min\bigl([a]\setminus(G_1\cup\cdots\cup G_{j-1})\bigr).
\end{equation}
In particular, $i_1=1$. In the prescribed order, examine the types $G$
that have not yet been tested and satisfy $m_G>0$ and $i_j\in G$.
If $L$ numbers of color $j$ have already been assigned when we test $G$,
match all $m_G$ tickets of this type, in numerical order, with the tickets
\[
 (j,L+1),\ \ldots,\ (j,L+m_G).
\]
If $G$ does not occur in $z$, keep these assignments and continue. Stop
at the first type $G$ that occurs in $z$ and set $G_j=G$.

\smallskip
\noindent\textit{After $q$ steps.}
List the unassigned tickets in the first bag in numerical order.
List the unassigned tickets in the second bag first by color, then
by number within each color. Match tickets in corresponding positions
in the two lists. Denote the resulting map by $\pi_z$.
The output sequence is
\[
 w=(\pi_z(z_1),\ldots,\pi_z(z_k)).
\]

\subsubsection{Unique recovery and counting sequences}

The bijection $\pi_z$ between the tickets in the two bags depends on
the covering sequence $z$. To compare the numbers of sequences, we show
that the output $w$ uniquely determines $\pi_z$ and hence $z$.

\begin{lemma}[Validity of the construction and unique recovery]\label{sdiv:lem:recover}
For each covering sequence $z$, all $q$ steps can be carried out. All
$q$ colors occur in the output sequence $w$. The sequence $w$ uniquely
determines the input $z$, including the ticket number at each position.
\end{lemma}

\begin{proof}
Divide the ordered probabilities $p_1,\ldots,p_a$ into $q$ groups of
$s$ consecutive entries. The first entry in each group is at most the
group average. Hence, by~\eqref{sdiv:eq:marginal} and~\eqref{sdiv:eq:colors},
\[
 \sum_{j=1}^q p_{(j-1)s+1}\le\frac1s\sum_{i=1}^a p_i=1,
 \qquad M_j\ge Np_{(j-1)s+1}\quad(1\le j\le q).
\]
At the beginning of step $j$, the selected types cover at most
$(j-1)s<a$ elements. Thus $i_j$ exists and $i_j\le(j-1)s+1$.
There are enough tickets of color $j$ even for all types containing
$i_j$, since
\[
 \sum_{G\ni i_j}m_G=Np_{i_j}\le Np_{(j-1)s+1}\le M_j .
\]
Since $z$ covers $[a]$, it contains a type that includes $i_j$.
This type cannot have been selected earlier, since then $i_j$ would
already be covered. Nor can it have been passed over, since only types
absent from $z$ are passed over. Thus $G_j$ will be found and at least
one position will receive color $j$. After $q$ steps, the two bags have
the same number of unassigned tickets. Matching them in the prescribed
order completes the bijection $\pi_z$.

We now recover this bijection from $w$. At the first step, $i_1=1$ is
known. Once the preceding steps have been recovered, the types already
tested are known, as is $i_j$ by~\eqref{sdiv:eq:pivot}. In the same
order, examine the untested types with $i_j\in G$ and $m_G>0$.
Starting at $1$, assign the next $m_G$ numbers of color $j$ to each
successive type $G$.

To recover $G_j$, find the first interval that contains the number of
a ticket of color $j$ occurring in $w$. The type assigned to that interval is $G_j$. Indeed, $z$ contains no tickets
of the types passed over, and the numbers assigned to those types
remain reserved for them. For each color, the final matching of the
unassigned tickets uses only numbers beyond the initial interval already
assigned. It therefore cannot affect the first interval encountered in
this search. At least one ticket of type $G_j$ occurs in $z$, so a
number from its interval occurs in $w$. This determines the stopping
point and all assignments at step $j$. After $q$ steps, the remaining
tickets are matched in the prescribed order. We now have
$z_t=\pi_z^{-1}(w_t)$ at every position $t$.
\end{proof}

Each bag gives $N^k$ equally likely sequences. By
Lemma~\ref{sdiv:lem:recover}, distinct covering sequences have distinct
images, so
\begin{equation}\label{sdiv:eq:injection}
 \frac{\#\{\text{covering sequences}\}}{N^k}
 \le
 \frac{\#\{\text{sequences containing all }q\text{ colors}\}}{N^k}.
\end{equation}

\paragraph{An example of encoding and recovery.}
Let $a=4$, $s=2$, and suppose that the four types $12,13,24,34$ are
equally likely, where $ij$ denotes the pair $\{i,j\}$. Take $N=4$ and
number the tickets in the order shown. All $p_i=1/2$, so the second bag
contains two tickets of each of two colors. Consider the covering
sequence $z=(2,3,4)$ with types $(13,24,34)$. At the first step,
$i_1=1$. Type $12$ is absent from the sequence, and the next type,
$13$, is present. At the second step, $i_2=2$, and $24$ is the first
eligible untested type. The last ticket is matched with the remaining
ticket of the second color, giving
\[
\begin{array}{c|cccc}
\text{ticket number}&1&2&3&4\\
\text{type}&12&13&24&34\\
\pi_z& (1,1)&(1,2)&(2,1)&(2,2)
\end{array}
\qquad
w=\bigl((1,2),(2,1),(2,2)\bigr).
\]
Conversely, $w$ contains $(1,2)$ but not $(1,1)$, so the first step
passed over $12$ and stopped at $13$. The next uncovered element is
therefore $2$. The presence of $(2,1)$ shows that the second step
stopped at $24$. Number $2$ of the second color was assigned to the
remaining ticket of type $34$ after type $24$ was chosen. Inverting
the resulting bijection recovers the original sequence $z=(2,3,4)$.

\subsubsection{Extension to arbitrary probabilities}

Now let the $y_G$ be arbitrary. As in~\eqref{sdiv:eq:colors}, define
the color probabilities by
\begin{equation}\label{sdiv:eq:shares}
 c_j=p_{(j-1)s+1}\quad(1\le j<q),\qquad
 c_q=1-\sum_{j=1}^{q-1}c_j .
\end{equation}
The bound on the sum of the first entries in the groups gives
$c_q\ge p_{(q-1)s+1}$. In the first bag, assign exactly
$m_G=\lfloor Ny_G\rfloor$ tickets to type $G$. Leave the remaining
tickets blank. A sequence counts as covering only if every ticket bears a type
and these types have union $[a]$. In the second bag, set $M_j=\lfloor Nc_j\rfloor$
for $j<q$ and $M_q=N-\sum_{j<q}M_j$. In particular, $M_q\ge Nc_q$.

The bound on the number of tickets required at step $j$ still holds
after rounding, since
\[
 \sum_{G\ni i_j}\lfloor Ny_G\rfloor
 \le\lfloor Np_{i_j}\rfloor\le\lfloor Nc_j\rfloor\le M_j .
\]
Blank tickets enter only in the final matching of unassigned
tickets. Thus the proof of Lemma~\ref{sdiv:lem:recover} still gives
\eqref{sdiv:eq:injection}. As $N\to\infty$, the left-hand side of this
inequality satisfies
\[
 \sum_{G_1\cup\cdots\cup G_k=[a]}
 \prod_{t=1}^k\frac{\lfloor Ny_{G_t}\rfloor}{N}
 \ \longrightarrow\ \P(T_y\le k),
\]
since the sum over type sequences is finite. The right-hand side is a
polynomial in the proportions $M_j/N$, and $M_j/N\to c_j$.
Let $J_1,\ldots,J_k$ be independent colors with $\P(J_t=j)=c_j$.
In the limit we obtain
\begin{equation}\label{sdiv:eq:comparison}
 \P(T_y\le k)\le\P\bigl(\{J_1,\ldots,J_k\}=[q]\bigr).
\end{equation}

Now partition $[a]$ into $q$ blocks $B_1,\ldots,B_q$ of size $s$ and
draw block $B_j$ with probability $c_j$. The collection is complete by
draw $k$ precisely when each block has been drawn at least once. Its
probability is the right-hand side of~\eqref{sdiv:eq:comparison}.
Since the probabilities $c_j$ in~\eqref{sdiv:eq:shares} depend only
on $y$, the same distribution on the blocks gives this bound for all
$k$. For $k<q$, both probabilities are zero.

\subsubsection{Equalizing the color probabilities}

\begin{lemma}[Collecting all colors]\label{sdiv:lem:balance}
For $k\ge q$, the probability of collecting all $q$ colors in $k$
independent draws is maximized precisely when $c_1=\cdots=c_q=1/q$.
\end{lemma}

\begin{proof}
The probability of complete collection is continuous on the probability
simplex and attains a maximum. It is positive when all colors are
equally likely and zero if any color has zero probability. Hence all
color probabilities are positive at a maximum.

Choose two colors with probabilities $u,v$. Fix the positions occupied
by these two colors and the specific colors in all other positions.
Suppose $h$ positions are assigned to the two colors. If $h\ge2$ and
all the other colors are present, the factor accounting for the
occurrence of both selected colors is $(u+v)^h-u^h-v^h$. The two
subtracted terms correspond to filling all $h$ positions with a single
color. When $h<2$ or one of the other colors is missing, the contribution
is zero.

If $u\ne v$, replacing both probabilities by their average strictly
increases every positive contribution, since
\[
 u^h+v^h>2\left(\frac{u+v}{2}\right)^h\qquad(h\ge2).
\]
For $k\ge q$, we can assign one position to each of the other $q-2$
colors, leaving $k-q+2\ge2$ positions for the two selected colors.
Since all color probabilities are positive, this gives a positive
contribution that increases strictly if $u\ne v$. Hence all color
probabilities are equal at a maximum.
\end{proof}

Applying Lemma~\ref{sdiv:lem:balance} to~\eqref{sdiv:eq:comparison} gives
\[
 \P(T_y\le k)
 \le \P\bigl(\{J_1,\ldots,J_k\}=[q]\bigr)
 \le \sum_{j=0}^q(-1)^j\binom qj\left(1-\frac{j}{q}\right)^k .
\]
The uniform distribution on the blocks attains the right-hand side,
which proves~\eqref{sdiv:eq:maximum}.

\subsection{Optimal distributions and cases of equality}

Suppose a distribution attains the maximum for a fixed $k\ge q$.
Then equality holds in~\eqref{sdiv:eq:comparison}, and all color
probabilities are $1/q$ by Lemma~\ref{sdiv:lem:balance}. In particular,
\eqref{sdiv:eq:shares} gives $p_1=c_1=1/q$. By~\eqref{sdiv:eq:marginal}
and our ordering of the elements,
\[
 \min_{i\in[a]}p_i=p_1=\frac1q=\frac{s}{a}
            =\frac1a\sum_{i=1}^a p_i .
\]
Since the minimum equals the average, $p_i=1/q$ for every $i$.

Suppose two distinct types $G,H$ of positive probability intersect.
We construct sequences containing all $q$ colors that have no covering
preimage. Since $|G|=|H|$, both $G\setminus H$ and $H\setminus G$ are
nonempty. As all $p_i$ are equal, we can renumber the elements so that
$1\in G\setminus H$ and $2\in H\setminus G$. Order the types with $G$
first and $H$ second. The bound~\eqref{sdiv:eq:injection} holds for
this order as well.

Take $N$ to be a sufficiently large multiple of $q$. Then $M_j=N/q$
for each $j$. At the first step, the first $m_G=\lfloor Ny_G\rfloor$
numbers of color $1$ are assigned to type $G$. If $G$ is chosen, then
$i_2=2$, and the first $m_H=\lfloor Ny_H\rfloor$ numbers of color $2$
are assigned to type $H$.

Form sequences from the second bag as follows. In each of the first
$k-q+1$ positions, choose any of the first $m_G$ tickets of color~$1$.
In the next position, choose any of the first $m_H$ tickets of color~$2$.
The last $q-2$ positions contain
tickets of colors $3,\ldots,q$, respectively, with arbitrary numbers.
Every such sequence contains all colors. When recovering a preimage,
the first interval of color $1$ would select type $G$, and the first
interval of color $2$ would select type $H$. Thus the first $k-q+1$
tickets of the preimage would have type $G$ and the next ticket would
have type $H$. Even with the last $q-2$ sets, the number of elements
covered would be at most
\[
 |G\cup H|+(q-2)s=qs-|G\cap H|<a .
\]
Hence none of these sequences has a covering preimage. Their proportion
among all $N^k$ sequences is
\[
 \left(\frac{m_G}{N}\right)^{k-q+1}
 \frac{m_H}{N}\left(\frac1q\right)^{q-2}
 \ \longrightarrow\ y_G^{\,k-q+1}y_H\left(\frac1q\right)^{q-2}>0.
\]
The right-hand side of~\eqref{sdiv:eq:injection} exceeds the left-hand
side by at least this proportion. Letting $N$ tend to infinity through
multiples of $q$ gives strict inequality in~\eqref{sdiv:eq:comparison},
a contradiction.

Thus the sets of positive probability are pairwise disjoint. Since
$p_i>0$ for every $i$, they cover $[a]$ and form a partition into $q$
blocks of size $s$. For every block $G$ and every $i\in G$, we have
$y_G=p_i=1/q$. This proves the equality statement in
Theorem~\ref{sdiv:thm:main}.

\subsection{Expected collection time}

\begin{theorem}[Minimum expected collection time when $s\mid a$]\label{sdiv:cor:mean}
Let $a=qs$, where $q\ge2$ is an integer. Then
\begin{equation}\label{sdiv:eq:mean}
 M(qs,(q-1)s)=qH_q,\qquad H_q=\sum_{j=1}^q\frac1j .
\end{equation}
The minimum is attained precisely when $D_t$ is chosen uniformly from
the blocks of a partition of $[a]$ into $q$ sets of size $s$.
\end{theorem}

\begin{proof}
Let $T_*$ be the collection time when blocks are chosen uniformly.
By Theorem~\ref{sdiv:thm:main}, for every $k\ge0$ we have
$\P(T_y>k)\ge\P(T_*>k)$. The tail-sum formula gives
\begin{equation}\label{sdiv:eq:tail}
 \E T_y=\sum_{k=0}^{\infty}\P(T_y>k)
 \ \ge\ \sum_{k=0}^{\infty}\P(T_*>k)=\E T_* .
\end{equation}
The formula also applies when the expected collection time is infinite.

Once $j$ blocks have been collected, the probability of obtaining a
new block on the next draw is $(q-j)/q$. The expected number of draws
until a new block appears is therefore $q/(q-j)$. Summing over $j$ gives
\[
 \E T_*=\sum_{j=0}^{q-1}\frac{q}{q-j}
              =q\sum_{r=1}^{q}\frac1r=qH_q .
\]

If $y$ is not uniform on the blocks of any partition of $[a]$ into
$s$-sets, Theorem~\ref{sdiv:thm:main} gives a strict tail inequality at
$k=q$. Hence the inequality in~\eqref{sdiv:eq:tail} is also strict.
\end{proof}

\subsection{Components of the support of the drawn sets}

Let $1\le s<a$. If a union of components has positive mass and contains
$ms$ vertices for some $m\ge2$, Theorem~\ref{sdiv:thm:main} applies to
the conditional distribution of the drawn sets on this union.

\begin{samepage}
\begin{theorem}[Components of the support of the drawn sets]\label{sdiv:cor:components}
Suppose a distribution attains $\alpha_{\max}(a,b,k)$ for some
$k\ge\lceil a/s\rceil$ or minimizes $\E T_y$. For any distinct
components $C_1,\ldots,C_\ell$ of the support, where $\ell\ge1$,
\[
 s\mid\bigl(|C_1|+\cdots+|C_\ell|\bigr)
 \quad\Longleftrightarrow\quad |C_1|=\cdots=|C_\ell|=s .
\]
In particular, if $s\mid a$, every component has size $s$.
There are at most $s-\gcd(a,s)$ components whose sizes are not divisible
by $s$. The sum of the sizes of any nonempty subfamily of these
components is also not divisible by $s$.
\end{theorem}
\end{samepage}

\begin{proof}
Let $r=\lceil a/s\rceil$. Partition $[a]$ into $r$ groups of size $s$,
except possibly for a smaller last group. If the last group has fewer
than $s$ elements, enlarge it to size $s$ using elements from the other
groups. Under the uniform distribution on the resulting $r$ sets,
there is a positive probability that all of them appear in $r$ draws.
Hence, for every $k\ge r$, the probability of complete collection under
a maximizing distribution is positive. At a minimizer of the mean
collection time, the mean is finite by Proposition~\ref{pre:tail}.
In either case every vertex belongs to an edge of the support.

Let $U$ be the union of the selected components and suppose that $|U|=ms$
for some integer $m\ge1$. We show that the edges in $U$ partition it
into blocks of size $s$. If $m=1$, any edge in $U$ equals $U$. Suppose
$m\ge2$, and write $W=[a]\setminus U$. Every set of positive probability
lies entirely in $U$ or entirely in $W$.

For $S\subseteq[a]$, write $P_h(S)$ for the probability that the union
of $h$ drawn sets $D_t$ is exactly $S$. In particular,
$P_0(\varnothing)=1$ and $P_0(S)=0$ for $S\ne\varnothing$. For every
$k\ge0$,
\begin{equation}\label{sdiv:eq:positions}
 P_k([a])=\sum_{h=0}^{k}\binom{k}{h}P_h(U)P_{k-h}(W).
\end{equation}
Fix $h$ of the $k$ positions. By independence of the draws, the
probability that the sets in these positions have union $U$ and those
in the remaining positions have union $W$ equals $P_h(U)P_{k-h}(W)$.
Each covering sequence determines the positions of its sets from $U$,
so it is counted exactly once. If $W=\varnothing$, only the term with
$h=k$ remains.

Keep the probabilities of all sets outside $U$ fixed and distribute
the total probability $\sigma>0$ of the sets in $U$ equally among the
$m$ blocks of a partition of $U$ into sets of size $s$. The factors
$P_{k-h}(W)$ in~\eqref{sdiv:eq:positions} stay the same. For $h<m$,
we have $P_h(U)=0$ under either distribution. If the union of $h$ sets
is $U$, all of them lie in $U$. The probability that all $h$ sets lie
in $U$ is $\sigma^h$. Conditional on this event they are independent,
and Theorem~\ref{sdiv:thm:main} applied to their common distribution on
$ms$ vertices gives
\[
 P_h(U)\le \sigma^h\,\alpha_{\max}\bigl(ms,(m-1)s,h\bigr)
 \qquad(h\ge m).
\]
The new distribution attains the right-hand side, so no term
in~\eqref{sdiv:eq:positions} decreases.

Suppose the original edges in $U$ did not form a partition. By the
equality condition in Theorem~\ref{sdiv:thm:main}, the inequality for
$P_h(U)$ was then strict for every $h\ge m$. Every term
in~\eqref{sdiv:eq:positions} that was positive before the replacement
therefore strictly increases. If the original distribution maximizes $P_k([a])$, then $P_k([a])>0$,
so the replacement strictly increases the probability of collecting all
coupons within $k$ draws. If the original distribution minimizes the
expected collection time, then $\E T_y<\infty$, so $P_k([a])>0$
for at least one $k$. The corresponding tail probability strictly
decreases and none of the others increases. By~\eqref{sdiv:eq:tail},
the expected collection time decreases as well. In either case we have a
contradiction. Thus the components in $U$ are blocks of size $s$.
The converse follows by adding their sizes.

Consequently, among the components whose sizes are not divisible by
$s$, no nonempty subfamily has total size divisible by $s$. The sum of
the sizes of all these components is congruent to $a$ modulo $s$.
Applying the counting argument at the end of the proof of
Theorem~\ref{bdiv:components}, with $s$ in place of $b$, gives the
bound $s-\gcd(a,s)$.
\end{proof}

\section{Missing sets when \texorpdfstring{$b=2$}{b=2}}\label{b2:section}

In this section the superscript $-$ refers to missing pairs. In the next
section the superscript $+$ refers to drawn pairs.

Let $a\ge3$ and $s=a-2$. We work with the missing pairs $A_t$.
For even $a=2q$, we call the uniform distribution on $q$ pairwise
disjoint pairs \emph{canonical}. For odd $a=2m+3$, consider a triangle
and $m$ isolated edges, all vertex-disjoint. We call a distribution on
this graph \emph{canonical} if the three edges of the triangle have
equal probabilities and the $m$ isolated edges have equal probabilities.
When $m=0$, this is the uniform distribution on the edges of a triangle.

Let $m\ge1$ and let $\theta\in[0,1]$ be the total probability assigned
to the edges of the triangle. Each of these edges has weight
$v=\theta/3$, and each isolated edge has weight $u=(1-\theta)/m$.
For $k\ge1$, the probability of collecting all coupons within $k$ draws
under this distribution is
\begin{equation}\label{b2:eq:canonical}
 F^{-}_{m,k}(\theta)=1-mu^k-3(2^k-1)v^k,\qquad
 u=\frac{1-\theta}{m},\quad v=\frac{\theta}{3}.
\end{equation}
The support is $K_3\sqcup mK_2$ when $0<\theta<1$.
Formula~\eqref{b2:eq:canonical} also covers the endpoints.
When $\theta=0$, only the $m$ isolated edges have positive weight,
and when $\theta=1$, only the edges of the triangle do.
For $k=0$, we set $F^{-}_{m,0}=0$.
Below we derive~\eqref{b2:eq:canonical} from the probability that the
missing pairs have nonempty intersection.

\begin{theorem}[Maximum probability of complete collection]\label{b2:thm:main}
If $a=2q$, $q\ge2$, then for every $k\ge2$,
\begin{equation}\label{b2:eq:even}
 \alpha_{\max}(2q,2,k)=1-q^{1-k}.
\end{equation}
If $a=2m+3$, $m\ge1$, then for every $k\ge2$,
\begin{equation}\label{b2:eq:odd}
 \begin{aligned}
 \alpha_{\max}(2m+3,2,k)
 &=\max_{0\le\theta\le1}F^{-}_{m,k}(\theta)
 =1-\frac{2^k-1}{(m\rho_k+3)^{k-1}},\\
 \rho_k&=(2^k-1)^{1/(k-1)}.
 \end{aligned}
\end{equation}
In this odd case, a canonical distribution maximizes the probability of
collecting all coupons within $k$ draws if and only if the total mass of
the triangle edges is
\begin{equation}\label{b2:eq:optimal-parameter}
 \theta_{m,k}:=\frac3{m\rho_k+3}.
\end{equation}
For $k\ge3$, in both cases every maximizer is obtained from the
specified canonical distribution by relabeling the coupons. In the odd
case, each isolated edge has weight $\rho_k/(m\rho_k+3)$ and each
triangle edge has weight $1/(m\rho_k+3)$.

For $a=3$, we have
$\alpha_{\max}(3,2,k)=1-(2^k-1)/3^{k-1}$ for $k\ge2$.
For $k\ge3$, only the uniform distribution on the edges of the triangle
attains the maximum. For $k=0,1$, complete collection is impossible for
every $a\ge3$.
\end{theorem}

Maximizing $\Phi_k$ for $k\ge3$ and minimizing the expected collection
time both reduce to minimizing functionals of the same form. The same
perturbations of the edge weights therefore characterize the supports of
the optimal distributions in both problems. For $k=2$,
Theorem~\ref{k2:main} gives all cases of equality.

\subsection{Structure of the optimal distributions}

\subsubsection{Probability of a nonempty intersection and expected collection time}

Let $x_{ij}$ be the probability that $\{i,j\}$ is the missing pair,
and let $p_i=\sum_{j\ne i}x_{ij}$ be the probability that coupon $i$
is missing from a single draw. Then
$\sum_{i<j}x_{ij}=1$ and $\sum_i p_i=2$.
For $k\ge1$, inclusion--exclusion applied to the events
$\{i\in A_1\cap\cdots\cap A_k\}$ gives
\begin{equation}\label{b2:eq:tail-probability}
 1-\Phi_k(x)
 =\P(T_x>k)
 =\sum_{i=1}^{a}p_i^k-\sum_{i<j}x_{ij}^k.
\end{equation}
A pair cannot contain three distinct coupons, so all higher-order terms
in inclusion--exclusion vanish.

If $p_i=1$ for some $i$, then coupon $i$ is never drawn and
$\E T_x=\infty$. If $p_i<1$ for all $i$, then
$\P(T_x>k)\le\sum_i p_i^k$, so the expected collection time is finite.
Summing~\eqref{b2:eq:tail-probability} gives
\begin{equation}\label{b2:eq:mean-functional}
 \begin{aligned}
 \E T_x
 &=1+\sum_i\frac{p_i}{1-p_i}
           -\sum_{i<j}\frac{x_{ij}}{1-x_{ij}}\\
 &=2+\sum_i\psi(p_i)-\sum_{i<j}\psi(x_{ij}),
 \qquad \psi(t)=\frac{t^2}{1-t}.
 \end{aligned}
\end{equation}
The second line uses $t/(1-t)=t+\psi(t)$ and
$\sum_i p_i-\sum_{i<j}x_{ij}=1$.

By Proposition~\ref{pre:tail}, both extrema are attained.
The minimum expected collection time is finite, so every minimizer
satisfies $p_i<1$ for all $i$.

We shall minimize
\begin{equation}\label{b2:eq:objective}
 \mathcal F(x)=\sum_i\varphi(p_i)-\sum_{i<j}\varphi(x_{ij}),
 \qquad
 \varphi(t)=t^k\ (k\ge3)\quad\text{or}\quad
 \varphi(t)=\psi(t).
\end{equation}
The first choice gives $\mathcal F(x)=1-\Phi_k(x)$, and the second gives
$\mathcal F(x)=\E T_x-2$ for distributions with finite expected
collection time. In both cases $\varphi(0)=\varphi'(0)=0$, and
$\varphi^{\prime\prime}$ is positive for $t>0$ and strictly increasing.
For the expected collection time we use
\[
 \psi(t)=\sum_{r=2}^{\infty}t^r,\qquad
 \psi^{\prime\prime}(t)=\frac2{(1-t)^3}\quad(0\le t<1).
\]

\subsubsection{Perturbing the weights and determining the possible components}

The support graph is the graph on $[a]$ whose edges are the pairs $ij$
with $x_{ij}>0$. The mass of a component $C$ is
$\sigma(C)=\sum_{e\in E(C)}x_e$.
The functional~\eqref{b2:eq:objective} is the sum of the contributions
from the components, and isolated vertices contribute zero.
In perturbing the weights we may add to the support a pair whose initial
weight is zero. Transferring mass from one edge to another will rule out
all components with leaves except isolated edges. Under a perturbation
that preserves all the $p_i$, the first sum in~\eqref{b2:eq:objective}
is constant. Strict convexity of the second sum will allow us to rule
out extra edges in components without leaves.

\begin{lemma}[Perturbations preserving the marginal probabilities]
\label{b2:lem:kernel}
Suppose that the edges of the support carry numbers $\delta_e$, not all
zero, such that $\sum_{e\ni i}\delta_e=0$ at every vertex $i$.
Then the distribution is not a local minimizer of~\eqref{b2:eq:objective}.
\end{lemma}

\begin{proof}
For small $|t|$, set $x_e(t)=x_e+t\delta_e$.
The probabilities $p_i$ are preserved, as is the total mass, since
$2\sum_e\delta_e=\sum_i\sum_{e\ni i}\delta_e=0$.
Along this perturbation,
\[
 \frac{d^2}{dt^2}\mathcal F(x+t\delta)
 =-\sum_e\delta_e^2\varphi^{\prime\prime}(x_e+t\delta_e)<0.
\]
Thus $t=0$ cannot be a local minimum.
\end{proof}

\begin{lemma}[Components of the support]\label{b2:lem:components}
At a local minimum of~\eqref{b2:eq:objective}, every component of the
support containing an edge is either an isolated edge or a simple odd cycle.
\end{lemma}

\begin{proof}
Suppose that $u$ is a leaf, $x_{uv}=c>0$, and $v$ has another neighbor
$w$ with $x_{vw}=d>0$. Then $p_u=c$ and $x_{uw}=0$.
Transfer mass $t$ from edge $vw$ to edge $uw$, where $0<t<d$.
This raises $p_u$ by $t$ and lowers $p_v$ by $t$, while $p_w$ stays
constant because both edges contain $w$. Denote the resulting
distribution by $x(t)$. Since $\varphi'(0)=0$ and $p_v\ge c+d$,
\[
 \left.\frac{d}{dt}\mathcal F(x(t))\right|_{t=0+}
 =\varphi'(c)-\varphi'(p_v)+\varphi'(d)
 \le \varphi'(c)+\varphi'(d)-\varphi'(c+d)<0.
\]
The last inequality is strict because $\varphi^{\prime\prime}$ is strictly
increasing and
\[
 \varphi'(c+d)-\varphi'(c)
 =\int_0^d\varphi^{\prime\prime}(c+r)\,dr
 >\int_0^d\varphi^{\prime\prime}(r)\,dr=\varphi'(d).
\]
For small $t>0$ the functional decreases, contradicting local
minimality. Thus a component with a leaf consists of a single edge.

In a component containing an edge but no leaves, all degrees are at
least two. If some degree is greater than two, then $|E(C)|>|V(C)|$.
The system $\sum_{e\ni i}\delta_e=0$ has more unknowns than equations
and hence has a nonzero solution, contrary to Lemma~\ref{b2:lem:kernel}.
Thus all degrees are two and the component is a cycle. On an even
cycle, assigning alternating values $+1,-1$ to the edges again yields
a perturbation excluded by the lemma. Only odd cycles remain.
\end{proof}

\subsubsection{Ruling out long cycles}

\begin{lemma}[A switch on a cycle]\label{b2:lem:long-cycle}
The support graph of a global minimizer of~\eqref{b2:eq:objective}
consists of isolated edges, triangles and possibly isolated vertices.
\end{lemma}

\begin{proof}
By Lemma~\ref{b2:lem:components}, it suffices to rule out cycles of
length at least five. Label four consecutive vertices of such a cycle so that the middle
edge has maximum weight $x_{23}=w$. Let $x_{12}=c$ and $x_{34}=d$.
Choosing the direction around the cycle so that $c\le d$, we have
$0<c\le d\le w$. Edge $14$ is absent.
Set
\begin{equation}\label{b2:eq:cycle-switch}
 x'_{12}=0,\qquad x'_{23}=w+c,\qquad
 x'_{34}=d-c,\qquad x'_{14}=c.
\end{equation}
All weights are nonnegative, and both the total mass and all the $p_i$
are preserved. Hence
\[
 \begin{aligned}
 \mathcal F(x')-\mathcal F(x)
 &=\varphi(w)+\varphi(d)-\varphi(w+c)-\varphi(d-c)\\
 &=-\int_0^c\bigl(\varphi'(w+t)-\varphi'(d-c+t)\bigr)\,dt<0.
 \end{aligned}
\]
Here $w>d-c$ and $\varphi'$ is strictly increasing. For the expected
collection time, preserving $p_i<1$ ensures that the new value is finite.
This contradiction proves the claim.
\end{proof}

\subsubsection{The number of components and their edge weights}

\begin{lemma}[Edge weights of a triangle]\label{b2:lem:triangle}
For a fixed positive total mass, a triangle makes the smallest
contribution to~\eqref{b2:eq:objective} precisely when its three edges
have equal weights. An isolated edge of mass $\sigma$ contributes
$\varphi(\sigma)$.
\end{lemma}

\begin{proof}
For an isolated edge of mass $\sigma$, the marginal probability $p_i$
equals $\sigma$ at each endpoint. Its contribution is therefore
$2\varphi(\sigma)-\varphi(\sigma)=\varphi(\sigma)$.

Let $w_1,w_2,w_3$ be the edge weights of a triangle with total mass
$\sigma$. At an interior minimum, transferring mass between $w_1$ and $w_2$
while keeping $w_3$ fixed gives zero first derivative. Canceling
the common term gives
\[
 \varphi'(w_1+w_3)-\varphi'(w_1)
 =\varphi'(w_2+w_3)-\varphi'(w_2).
\]
For $w_3>0$, the function $t\mapsto\varphi'(t+w_3)-\varphi'(t)$ is
strictly increasing because its derivative is
$\varphi^{\prime\prime}(t+w_3)-\varphi^{\prime\prime}(t)>0$.
Thus $w_1=w_2$, and by symmetry all three weights are equal.

On the boundary, say when $w_3=0$, the contribution is
$\varphi(\sigma)$. For $\varphi(t)=t^k$, $k\ge3$, equal weights give
\[
 3\left[\left(\frac{2\sigma}{3}\right)^k
             -\left(\frac{\sigma}{3}\right)^k\right]
 =\frac{2^k-1}{3^{k-1}}\,\sigma^k<\sigma^k.
\]
The inequality $2^k-1<3^{k-1}$ follows by induction starting from
$7<9$. For $\varphi=\psi$ and $\sigma<1$, sum these comparisons over
$k\ge2$, using equality for $k=2$ and strict inequality for $k\ge3$.
For $\varphi=\psi$ and $\sigma=1$, the expected collection time is
infinite on the boundary and finite for equal weights. Hence the minimum is attained
in the interior and is unique.
\end{proof}

\begin{lemma}[Canonical support]\label{b2:lem:canonical-support}
Suppose that $x$ maximizes $\Phi_k$ for some $k\ge3$ or minimizes
$\E T_x$. For even $a$, its support is a partition into pairs.
For odd $a$, it consists of a triangle and isolated edges on the
remaining vertices. In both cases all isolated edges are equally likely,
and the triangle edges, when present, have equal weights.
\end{lemma}

\begin{proof}
Both problems reduce to finding a global minimum
of~\eqref{b2:eq:objective}. By Lemma~\ref{b2:lem:long-cycle}, the only
possible components are isolated edges, triangles and isolated
vertices. By Theorem~\ref{bdiv:components} with $b=2$, the number of
components of odd order is at most $2-\gcd(a,2)$. Hence for even $a$ all
components are isolated edges. For odd $a$, exactly one component has
odd order and is either a triangle or an isolated vertex.

In the latter case, take the isolated vertex and an isolated edge of
mass $\sigma>0$ and form a triangle on their three vertices.
Replace the weights $(\sigma,0,0)$ by $(\sigma/3,\sigma/3,\sigma/3)$.
By Lemma~\ref{b2:lem:triangle}, the contribution of these vertices
to~\eqref{b2:eq:objective} strictly decreases for both choices of
$\varphi$. All other components make the same contribution as before,
contradicting minimality.

Thus for odd $a$ the unique component of odd order is a triangle.
By Lemma~\ref{b2:lem:triangle}, its edges have equal weights. For
a fixed total mass on the isolated edges, strict convexity of
$\varphi$ implies that their contribution $\sum_e\varphi(x_e)$ is
minimal only when their weights are equal.
\end{proof}

\subsection{Optimal weights and cases of equality}

By Lemma~\ref{b2:lem:canonical-support}, for odd $a\ge5$ it remains
to determine the total mass $\theta$ of the triangle.

\subsubsection{Probability of complete collection on the canonical support}

For $a=2q$, the canonical distribution is uniform on the $q$ pairs
of a fixed partition of $[a]$. The collection is incomplete after
$k\ge1$ draws if and only if the same pair is missing on all $k$ draws.
This event has probability $q^{1-k}$, so
Lemma~\ref{b2:lem:canonical-support} gives~\eqref{b2:eq:even}
for $k\ge3$.

Let $a=2m+3$, $m\ge1$. Under the canonical distribution $x(\theta)$
with total mass $\theta$ on the triangle, each of the $m$ isolated
edges has weight $u=(1-\theta)/m$ and each triangle edge has weight
$v=\theta/3$. By~\eqref{b2:eq:tail-probability}, the probability of
incomplete collection after $k\ge1$ draws is given by the right-hand
side of~\eqref{b2:eq:objective} with $\varphi(t)=t^k$.
An isolated edge contributes $2u^k-u^k=u^k$. Two edges of weight $v$
meet at each vertex of the triangle, so the triangle contributes
$3(2v)^k-3v^k$. Hence the probability of complete collection is
\[
 F^{-}_{m,k}(\theta)
 =1-\mathcal F(x(\theta))
 =1-mu^k-3(2^k-1)v^k,
\]
which proves~\eqref{b2:eq:canonical}.

For fixed $m$ and $k\ge2$, we determine the total mass $\theta$ on the
triangle that maximizes this probability.
Differentiating~\eqref{b2:eq:canonical} with respect to $\theta$ gives
\begin{equation}\label{b2:eq:derivatives}
 \frac{d}{d\theta}F^{-}_{m,k}(\theta)
 =k\bigl[u^{k-1}-(2^k-1)v^{k-1}\bigr].
\end{equation}
The right-hand side is strictly decreasing in $\theta$, so $F^{-}_{m,k}$
is strictly concave. The derivative vanishes when $u=\rho_k v$, where
$\rho_k=(2^k-1)^{1/(k-1)}$. Using $mu+3v=1$, we find that the unique maximum occurs at
\[
 u=\frac{\rho_k}{m\rho_k+3},\qquad
 v=\frac1{m\rho_k+3},\qquad
 \theta_{m,k}=3v=\frac3{m\rho_k+3}.
\]
At this point $u^{k-1}=(2^k-1)v^{k-1}$, so~\eqref{b2:eq:canonical}
gives
\[
 F^{-}_{m,k}(\theta_{m,k})
 =1-(mu+3v)u^{k-1}=1-u^{k-1}
 =1-\frac{2^k-1}{(m\rho_k+3)^{k-1}}.
\]
By Lemma~\ref{b2:lem:canonical-support}, this proves~\eqref{b2:eq:odd}
for $k\ge3$. When $a=3$ and $k\ge3$, the uniform distribution on the
edges of the triangle is the unique maximizer. Substituting $v=1/3$ in
the expression for its contribution gives
$\alpha_{\max}(3,2,k)=1-(2^k-1)/3^{k-1}$.

\subsubsection{Classification of the maximizers}

For $k\ge3$, all maximizers are canonical by
Lemma~\ref{b2:lem:canonical-support}. Thus for even $a$ the
distribution is uniform on a perfect matching, and for $a=3$ it is
uniform on the edges of a triangle. For odd $a\ge5$, the unique
optimal triangle mass is given by~\eqref{b2:eq:optimal-parameter}.

For $k=2$, Theorem~\ref{k2:main} describes all maximizers, and the
maximum is $1-1/\lfloor a/2\rfloor$. It is attained by choosing one
of $\lfloor a/2\rfloor$ disjoint pairs uniformly at random. For
$a=2m+3\ge5$, substituting $\rho_2=3$ in~\eqref{b2:eq:optimal-parameter}
gives $\theta_{m,2}=1/(m+1)$, and by~\eqref{b2:eq:canonical} the
canonical distribution with this mass also gives $1-1/(m+1)$.
For $a=3$, every distribution gives $\Phi_2=0$. This completes the
proof of Theorem~\ref{b2:thm:main}.

\subsection{Expected collection time}

\begin{theorem}[Minimum expected collection time]\label{b2:cor:mean}
For $a=2q$, $q\ge2$, we have
\begin{equation}\label{b2:eq:mean-even}
 M(2q,2)=2+\frac1{q-1}.
\end{equation}
The minimum is attained precisely when the missing pair is chosen
uniformly from the $q$ pairs of a partition of $[a]$. For $a=3$,
$M(3,2)=11/2$, and the uniform distribution on the three pairs is the
unique minimizer.

Let $a=2m+3$ with $m\ge1$. The expected collection time under the
canonical distribution $x(\theta)$ with total mass $\theta\in(0,1)$
on the triangle edges is
\begin{equation}\label{b2:eq:mean-odd}
 E_m^{-}(\theta)\coloneqq \E T_{x(\theta)}
 =1+\frac{mu}{1-u}+\frac{3v}{(1-2v)(1-v)},
 \qquad u=\frac{1-\theta}{m},\quad v=\frac{\theta}{3}.
\end{equation}
The minimum $M(2m+3,2)=E_m^{-}(\theta_{m,E})$ is attained precisely
by the canonical distributions with total triangle mass $\theta_{m,E}$.
Here $\theta_{m,E}$ is the unique $\theta\in(0,1)$ for which
\begin{equation}\label{b2:eq:mean-parameter}
 \frac1{(1-u)^2}=\frac2{(1-2v)^2}-\frac1{(1-v)^2}.
\end{equation}
Moreover,
\[
 \frac1{m+1}<\theta_{m,E}<\frac3{2m+3}.
\]
\end{theorem}

\begin{proof}
For $a=2q$, the assertion follows from Theorem~\ref{bdiv:mean} with
$b=2$. Now suppose that $a$ is odd. Proposition~\ref{pre:tail} ensures that the minimum is attained,
and Lemma~\ref{b2:lem:canonical-support} shows that every minimizer
is canonical. If $a=2m+3$ with $m\ge1$, the
support of a minimizer contains both the triangle and all $m$ isolated
edges, so $0<\theta<1$.

For the canonical distributions under consideration, $p_i=u<1$ at
the endpoints of the isolated edges and $p_i=2v<1$ on the triangle.
Thus with $\varphi=\psi$,
\eqref{b2:eq:mean-functional}--\eqref{b2:eq:objective} give
$\E T_x=2+\mathcal F(x)$, where $\psi(t)=t^2/(1-t)$.
An isolated edge contributes $2\psi(u)-\psi(u)=\psi(u)$ to
$\mathcal F(x)$, and the triangle contributes $3[\psi(2v)-\psi(v)]$.
For $a=3$ the minimizer has $v=1/3$, giving
$M(3,2)=2+3[\psi(2/3)-\psi(1/3)]=11/2$.
For $m\ge1$, summing over the components gives
\[
 E_m^{-}(\theta)=2+m\psi(u)+3[\psi(2v)-\psi(v)].
\]
Substituting $\psi(t)=t/(1-t)-t$ and using $mu+3v=1$
gives~\eqref{b2:eq:mean-odd}.

By the tail-sum formula~\eqref{pre:tail-equation},
\[
 E_m^{-}(\theta)=\sum_{k\ge0}\bigl(1-F^{-}_{m,k}(\theta)\bigr).
\]
The terms for $k=0,1$ are constant, and those for $k\ge2$ are strictly
convex by~\eqref{b2:eq:derivatives}. The sum is finite for each
$0<\theta<1$ by~\eqref{b2:eq:mean-odd}, so $E_m^{-}$ is strictly
convex on $(0,1)$. Differentiating the expression for $E_m^{-}$ in
terms of $\psi$ with respect to $\theta$, we obtain
\[
 \frac{d}{d\theta}E_m^{-}(\theta)=2\psi'(2v)-\psi'(v)-\psi'(u).
\]
At the minimum this derivative vanishes. Since
$\psi'(t)=(1-t)^{-2}-1$, we obtain~\eqref{b2:eq:mean-parameter}.
By strict convexity its solution $\theta_{m,E}$ is unique.
By~\eqref{b2:eq:mean-odd}, every canonical distribution with this mass
attains the same minimum.

The function $\psi'$ is strictly increasing and strictly convex on
$[0,1)$. Hence at the minimum,
\[
 \psi'(2v)<2\psi'(2v)-\psi'(v)=\psi'(u)<\psi'(3v).
\]
The first inequality follows from $v<2v$, and the last follows from
strict convexity at the midpoint of $[v,3v]$. Thus $2v<u<3v$.
Substituting $u=(1-\theta)/m$ and $v=\theta/3$ gives the stated bounds
on $\theta_{m,E}$.
\end{proof}

\section{Drawn sets when \texorpdfstring{$s=2$}{s=2}}\label{s2:sec:main}

Let $a\ge3$. The draws $D_1,D_2,\ldots$ are independent with a common
distribution $y=(y_e)$ on the two-element subsets of $[a]$. Set
\[
 T_y\coloneqq\min\{k:D_1\cup\cdots\cup D_k=[a]\},\qquad
 \Phi_k(y)\coloneqq\mathbb P(T_y\le k).
\]
The support graph of $y$ has vertex set $[a]$ and an edge for each pair
drawn with positive probability. We call a component consisting of a
single edge ($K_2$) an \emph{isolated edge}.
In the notation for missing sets, $b=a-2$ and
$y_D=x_{[a]\setminus D}$ for $D\in\binom{[a]}2$. Hence
\[
 \alpha_{\max}(a,a-2,k)=\max_y\Phi_k(y),\qquad
 M(a,a-2)=\min_y\mathbb E T_y.
\]
We set $\Phi_0=0$ and use the convention $0^0=1$ in the formulas.
All uniqueness statements are understood up to a permutation of the coupons.

For even $a=2q$, we call the uniform distribution on a perfect matching
canonical. For $a=3$, the canonical distribution is uniform on the edges of a triangle.
For $a=2m+3\ge5$, the canonical support consists of a triangle and $m$
isolated edges. We call a distribution on this support canonical if all
isolated edges have the same mass and all three edges of the triangle have
the same mass.
Let $\theta$ denote the total mass of the triangle. Each of its edges then
has weight $v$, and each isolated edge has weight $u$, where
\begin{equation}\label{s2:eq:weights}
 v=\frac{\theta}{3},\qquad u=\frac{1-\theta}{m},\qquad 0\le\theta\le1.
\end{equation}
For $0<\theta<1$, the support graph is $K_3\sqcup mK_2$. When $\theta=0$,
the three vertices of the triangle are isolated, and when $\theta=1$, the
$2m$ vertices outside the triangle are isolated. Write
$F_{m,k}^{+}(\theta)$ for the probability of collecting all coupons within
$k$ draws and $E_m^{+}(\theta)$ for the expected collection time. At both
endpoints the probability is zero and the expected collection time is
infinite. The superscript $+$ refers to the case $s=2$.
In Section~\ref{s2:sec:profiles}, we use inclusion--exclusion to obtain
\begin{equation}\label{s2:eq:canonical}
 \begin{aligned}
 F_{m,k}^{+}(\theta)=\sum_{j=0}^{m}(-1)^j\binom mj
 \biggl[&(1-ju)^k-3\left(1-ju-\frac{2\theta}{3}\right)^k\\
        &+2(1-ju-\theta)^k\biggr].
 \end{aligned}
\end{equation}
To compare distributions simultaneously for all $k$, we introduce the interval
\begin{equation}\label{s2:eq:interval}
 I_m^{+}=\left[\frac{3}{2m+3},\frac{2}{m+2}\right].
\end{equation}

\begin{theorem}[Simultaneous comparison]\label{s2:thm:dom}
For every distribution $y$ on the drawn pairs, there exists a canonical
distribution $y^*$ such that
\[
 \Phi_k(y)\le\Phi_k(y^*)\qquad\text{simultaneously for all }k\ge0.
\]
For odd $a=2m+3\ge5$, we may choose $y^*$ with the weights
in~\eqref{s2:eq:weights} and the same $\theta\in I_m^{+}$ for all $k$.
\end{theorem}

For each $y$, every tail probability for $T_{y^*}$ is at most the
corresponding tail probability for $T_y$. Thus Proposition~\ref{pre:transfer} gives
$\mathbb E T_{y^*}\le\mathbb E T_y$. In particular, some canonical
distribution minimizes the expected collection time.
Theorem~\ref{s2:thm:structure} shows that every minimizer is canonical.
For odd $a\ge5$, the triangle mass that maximizes $\Phi_k$ depends on the
number of draws (Theorem~\ref{s2:thm:horizon}).

\begin{samepage}
\begin{theorem}[Maximum probability]\label{s2:thm:horizon}
For $a=2q\ge4$ and every $k\ge0$,
\begin{equation}\label{s2:eq:evenvalue}
 \alpha_{\max}(2q,2q-2,k)
 =\sum_{j=0}^{q}(-1)^j\binom qj\left(1-\frac jq\right)^k.
\end{equation}
For $k\ge q$, the maximizers are precisely the uniform distributions on
perfect matchings. For $k<q$, the probability is zero for every distribution.

For $a=2m+3\ge5$,
\[
 \alpha_{\max}(2m+3,2m+1,k)=\max_{\theta\in I_m^{+}}F_{m,k}^{+}(\theta).
\]
If $k\ge m+2$, all maximizers are canonical and have the same triangle
mass $\theta_{m,k}^{+}$. This is the unique solution in $(0,1)$ of
$\frac{\mathrm d}{\mathrm d\theta}F_{m,k}^{+}(\theta)=0$, and
\[
 \frac{3}{2m+3}<\theta_{m,k}^{+}\le\frac{2}{m+2},
\]
with equality on the right exactly when $k=m+2$.
For $k<m+2$, every distribution is optimal.

For $a=3$, we have $\alpha_{\max}(3,1,0)=0$ and
$\alpha_{\max}(3,1,k)=1-3^{1-k}$ for $k\ge1$.
For $k\ge2$, the unique maximizer is the uniform distribution on the
edges of the triangle. For $k=0,1$, every distribution is optimal.
\end{theorem}
\end{samepage}

\begin{theorem}[Minimum expected collection time]\label{s2:thm:mean}
For $a=2q\ge4$,
\[
 M(2q,2q-2)=qH_q,\qquad H_q=\sum_{j=1}^{q}\frac{1}{j}.
\]
The minimizers are precisely the uniform distributions on perfect matchings.
For $a=3$, we have $M(3,1)=5/2$, and only the uniform distribution on
the edges of the triangle attains the minimum.

For $a=2m+3\ge5$, all minimizers are canonical. Under each minimizer,
the probability that the drawn edge belongs to the triangle is
$\theta_{m,E}^{+}$, where
\[
 \frac{3}{2m+3}<\theta_{m,E}^{+}<\frac{2}{m+2}.
\]
This is the unique minimizer of $E_m^{+}$ on $(0,1)$ and the unique
solution of $\frac{\mathrm d}{\mathrm d\theta}E_m^{+}(\theta)=0$ on this
interval. Thus
\[
 M(2m+3,2m+1)=E_m^{+}(\theta_{m,E}^{+}).
\]
An explicit expression for $E_m^{+}$ is given in~\eqref{s2:eq:mean-product}.
\end{theorem}

For even $a=2q$, we use Theorems~\ref{sdiv:thm:main}
and~\ref{sdiv:cor:mean}. For odd $a\ge5$, we first determine the optimal
support and then compare distributions by the total mass of the triangle.

\subsection{Local comparisons}
Write $N(v)$ for the open neighborhood of a vertex $v$. Whenever a new
edge coincides with an existing edge, their masses are added.

\begin{lemma}[Isolating an edge]\label{s2:lem:stars}
Let $P,Q$ be two vertices such that $L=N(P)\setminus\{Q\}$ and
$R=N(Q)\setminus\{P\}$ are nonempty and disjoint. Let $\xi_i=y_{Pi}$,
$\zeta_j=y_{Qj}$, and $w=y_{PQ}\ge0$. Denote the total masses of the two
stars by
\[
 X\coloneqq\sum_{i\in L}\xi_i,\qquad Z\coloneqq\sum_{j\in R}\zeta_j
\]
and interchange $P$ and $Q$ if necessary so that $X\le Z$.
Delete the edges $Pi,Qj$, add mass $X$ to the edge $PQ$, and add mass
$\xi_i\zeta_j/X$ to each edge $ij$ with $(i,j)\in L\times R$.
Under the resulting distribution $\widetilde y$, the edge $PQ$ has mass $w+X$
and forms a component on its own. For $k\ge1$,
\[
 \Phi_k(\widetilde y)-\Phi_k(y)\ge kX\,p_{k-1},\qquad
 p_j\coloneqq\P_y\!\left(\bigcup_{i=1}^jD_i=[a]\setminus\{P,Q\}\right).
\]
In particular, suppose that $[a]\setminus\{P,Q\}$ can be covered by $h-1$
original edges avoiding $P$ and $Q$, where $h=\lceil a/2\rceil$. Then the
probability of complete collection strictly increases for every $k\ge h$.
\end{lemma}

\begin{figure}[H]
\centering
\setlength{\unitlength}{1mm}
\begin{picture}(160,57)
\put(20,44){\circle*{1.5}}
\put(60,44){\circle*{1.5}}
\put(20,12){\circle*{1.5}}
\put(60,12){\circle*{1.5}}
\put(20,44){\line(1,0){40}}
\put(20,12){\line(0,1){32}}
\put(60,12){\line(0,1){32}}
\put(20,49){\makebox(0,0){$P$}}
\put(60,49){\makebox(0,0){$Q$}}
\put(40,48){\makebox(0,0){$w$}}
\put(13,28){\makebox(0,0){$\xi_i$}}
\put(67,28){\makebox(0,0){$\zeta_j$}}
\put(20,6){\makebox(0,0){$i\in L$}}
\put(60,6){\makebox(0,0){$j\in R$}}
\put(77,29){\vector(1,0){14}}
\put(105,44){\circle*{1.5}}
\put(145,44){\circle*{1.5}}
\put(105,12){\circle*{1.5}}
\put(145,12){\circle*{1.5}}
\put(105,44){\line(1,0){40}}
\put(105,12){\line(1,0){40}}
\put(105,49){\makebox(0,0){$P$}}
\put(145,49){\makebox(0,0){$Q$}}
\put(125,48){\makebox(0,0){$w+X$}}
\put(125,20){\makebox(0,0){$\xi_i\zeta_j/X$}}
\put(105,6){\makebox(0,0){$i\in L$}}
\put(145,6){\makebox(0,0){$j\in R$}}
\end{picture}
\caption{Isolating the edge $PQ$. For each $(i,j)\in L\times R$,
add $\xi_i\zeta_j/X$ to the weight of $ij$.}
\end{figure}
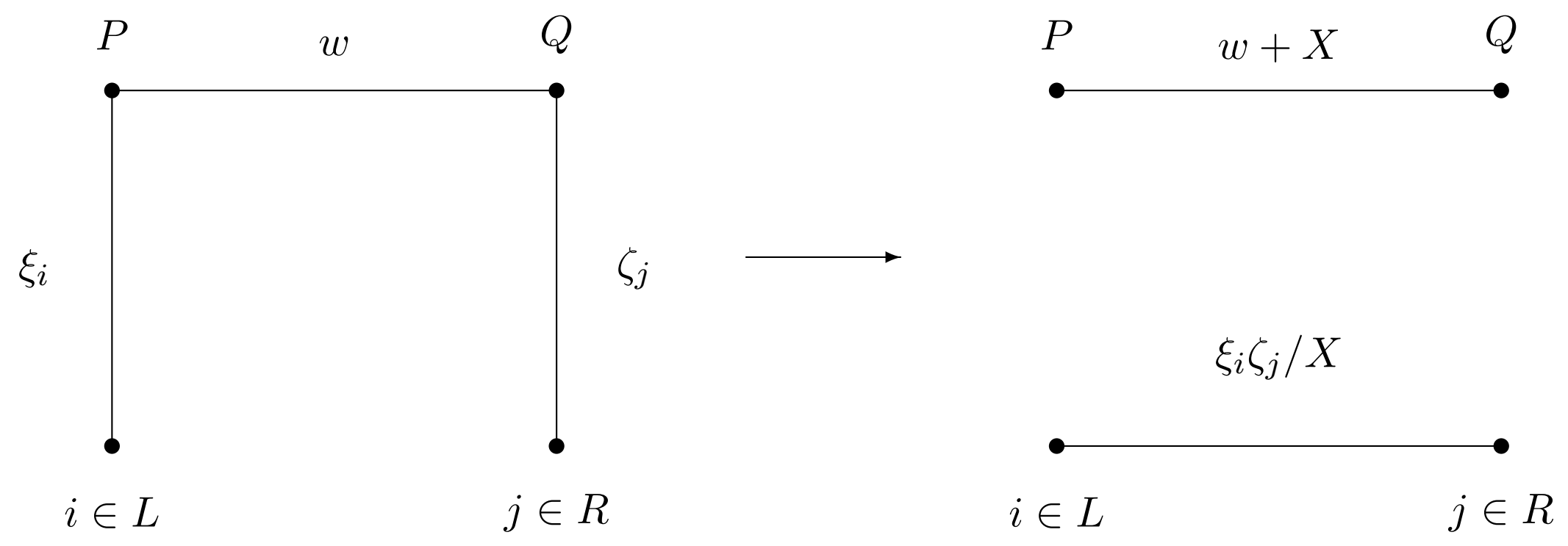

\begin{proof}
The total mass added to the edges between $L$ and $R$ is
$(\sum_i\xi_i)(\sum_j\zeta_j)/X=Z$, and the mass added to $PQ$ is $X$.
Condition on the positions and outcomes of all draws from the unchanged
part of the distribution, including the original mass $w$ of $PQ$. When
an old and a new outcome yield the same edge, we distinguish them by their
labels. Each such conditioning event has the same probability under both
distributions.

At each remaining position, choose a label shared by the two distributions.
Choose $X$ with probability $X/(X+Z)$ and $Z$ with probability $Z/(X+Z)$,
independently across positions.
A draw labeled $X$ changes from $Pi$ to $PQ$, and a draw labeled $Z$
changes from $Qj$ to $ij$. Let $N_X,N_Z$ be the numbers of occurrences
of these labels. The number of draws from $L$ changes from $N_X$
to $N_Z$. In both cases there are $N_Z$ draws from $R$.

Let $A_n$ be the probability of collecting all the remaining vertices of
$L$ in $n$ draws with probabilities $\xi_i/X$, and let $B_n$ be the
corresponding probability for $R$ with probabilities $\zeta_j/Z$.
The endpoints of a new edge $ij$ are chosen independently, since
$(\xi_i\zeta_j/X)/Z=(\xi_i/X)(\zeta_j/Z)$. Thus the probability
of covering both sides is the product of the two probabilities.

If any vertex outside $\{P,Q\}\cup L\cup R$ remains uncovered, both
distributions give zero probability of complete collection. If none of
the fixed draws used the original $PQ$ label of weight $w$ and one of
the counts is zero, complete
collection was impossible before the change, so the gain is nonnegative.
In all other cases, $P,Q$ are covered under both distributions.
When $N_X=N_Z$, the probabilities of complete collection are equal.

\begin{samepage}
Let $r<\ell$. Group together all arrangements of labels for which
$(N_X,N_Z)$ equals $(r,\ell)$ or $(\ell,r)$. The conditional probabilities
of complete collection in the two cases are
\[
 \begin{array}{c|ccc}
 (N_X,N_Z) & \text{before} & \text{after} & \text{gain} \\ \hline
 (r,\ell) & A_r B_\ell & A_\ell B_\ell & (A_\ell-A_r)B_\ell \\
 (\ell,r) & A_\ell B_r & A_r B_r & -(A_\ell-A_r)B_r
 \end{array}
\]
Both rows have the same number of label arrangements, so the ratio of
their probabilities is $\lambda=X^rZ^\ell/(X^\ell Z^r)=(Z/X)^{\ell-r}\ge1$.
Conditional on one of the two rows occurring, their probabilities are
$\lambda/(1+\lambda)$ and $1/(1+\lambda)$. The expected gain is
\[
 \begin{aligned}
 \Delta_{r,\ell}
 &=\frac{\lambda}{1+\lambda}(A_\ell B_\ell-A_r B_\ell)
   +\frac{1}{1+\lambda}(A_r B_r-A_\ell B_r)\\
 &=\frac{(A_\ell-A_r)(\lambda B_\ell-B_r)}{1+\lambda}\ge0.
 \end{aligned}
\]
Indeed, additional draws can only increase the probability of complete
collection, so $A_\ell\ge A_r$ and $\lambda B_\ell\ge B_\ell\ge B_r$.
\end{samepage}
Averaging over all cases gives $\Phi_k(\widetilde y)\ge\Phi_k(y)$.

For the quantitative bound, consider the conditioning events in which
exactly one position is free and the union of the other $k-1$ draws is
$[a]\setminus\{P,Q\}$. Complete collection is impossible before the
change. After the change it occurs precisely when the label is $X$,
giving $PQ$. For each possible position of this label, the contribution
to the gain is $X p_{k-1}$. These events are mutually exclusive, and the
other conditioning events give a nonnegative gain. Summing yields the
lower bound $kX p_{k-1}$. If the remaining vertices can be covered by
$h-1$ original edges avoiding $P$ and $Q$, then for $k\ge h$ we can
repeat edges of this cover to obtain a sequence of $k-1$ draws. Hence
$p_{k-1}>0$.
\end{proof}

\begin{lemma}[Eliminating a four-cycle]\label{s2:lem:cycle}
Let $C$ be a four-cycle in the support graph of $y$. Add $t,-t,t,-t$ to
the edge weights in cyclic order, leaving all other weights fixed. At
each endpoint of the feasible interval, at least one edge of the cycle
disappears. For each fixed $k$, there is an endpoint at which the
probability of complete collection is at least as large as under the
original distribution. More precisely, write $y_-,y_+$ for the
distributions at the endpoints, and let $y=(1-\lambda)y_-+\lambda y_+$,
where $0<\lambda<1$. Then
\[
 \Phi_k(y)\le (1-\lambda)\Phi_k(y_-)+\lambda\Phi_k(y_+)
 \qquad(k\ge0).
\]
If the support graph of at least one of $y_-,y_+$ consists of isolated
edges and at most one triangle that together cover all $a$ vertices, then
the inequality is strict for every $k\ge h=\lceil a/2\rceil$.
\end{lemma}

\begin{proof}
Condition on the positions of the draws from $C$ and on the edges drawn
at the remaining positions. The probability of this event is constant
in $t$, since the total mass of the cycle is preserved. If there are no
draws from $C$, or if some coupon outside the cycle remains uncollected,
the outcome is already determined. Suppose there are $r\ge1$ draws from
the cycle, and let $R$ be the set of its vertices still to be collected.
Write $p_e(t)$ for the conditional probability of drawing $e$ given that
the draw is from $C$. For each vertex $v$ of the cycle, the sum $b_v$
of the conditional probabilities of the two incident edges is constant
because their weights change by opposite amounts.

The probability of missing $v$ in all $r$ draws is $(1-b_v)^r$. To miss
both endpoints of an edge $e$, all $r$ draws must give the opposite edge
$\bar e$. It is impossible to miss two opposite vertices or any three
vertices. By inclusion--exclusion, the conditional probability of
collecting all of $R$ is
\[
 1-\sum_{v\in R}(1-b_v)^r
   +\sum_{\substack{e\in E(C)\\ e\subseteq R}} p_{\bar e}(t)^r .
\]
Each $p_{\bar e}(t)$ is affine and nonnegative, so its $r$th power is
convex. Averaging over the conditioning events shows that $\Phi_k(y(t))$
is convex, which gives the desired comparison with the endpoints.

Consider an endpoint whose support graph has the form stated in the
lemma. Replacing $t$ by $-t$ if necessary, take this to be the right
endpoint. Denote by $12$ and $34$ the opposite edges whose weights
increase with $t$. At the right endpoint they lie in distinct components,
and at least one, say $34$, is an isolated edge. Choose a cover of size
$h$ containing $12$ and $34$ by taking all isolated edges and, if there
is a triangle, two of its edges. Fix the first $h-2$ draws to be the
remaining edges of the cover in a prescribed order, and require the last
$r=k-h+2\ge2$ draws to come from the cycle. All vertices outside the cycle
are now covered, and neither $3$ nor $4$ has appeared. The only way to
miss both is to draw $12$ every time. Thus the inclusion--exclusion
formula contains the term $+p_{12}(t)^r$. Let $w_{12}$ be the original
weight of $12$, and let $W>0$ be the constant total mass of the cycle.
Then $p_{12}(t)=(w_{12}+t)/W$ is affine and strictly increasing in $t$.
For $r\ge2$, the function $x\mapsto x^r$ is strictly convex on
$[0,\infty)$, so $p_{12}(t)^r$ is strictly convex. The other terms are
constant or convex. The conditioning event has constant positive
probability, so the inequality remains strict after averaging.
\end{proof}

\subsection{Ruling out a bouquet of triangles}
\label{s2:sec:bouquet}

\begin{lemma}[Ruling out a bouquet]\label{s2:lem:bouquet}
Suppose that the support graph of a distribution on pairs is a friendship
graph, that is, a bouquet of $\nu\ge2$ triangles with a common vertex
and otherwise disjoint vertex sets. There is a distribution on pairs of the same $2\nu+1$
coupons whose support contains an isolated edge and for which the probability
of complete collection is strictly greater for every $k\ge\nu+1$.
For smaller $k$, both probabilities are zero.
\end{lemma}

\begin{proof}
Let $o$ denote the center of the bouquet. Call the other two vertices of
each petal its private vertices, the edge joining them its base, and the
two edges incident to $o$ its spokes. We distinguish two cases according
to the weights of these edges.

\paragraph{A light base.}
Suppose that the base of petal $oxy$ has weight $h$ no greater than the total
weight $w_S$ of the spokes of the other petals. Give each edge a label drawn
with the probability assigned to that edge. Denote the labels of
$xy,ox,oy$ by $H,X,Y$, respectively, and choose a private vertex $z$
of another petal. Keep the label probabilities fixed and change the pairs
they produce. Both $H$ and $X$ now produce $ox$, while $Y$ produces $yz$
and the label of each external spoke $ow$ produces $yw$. The external bases
retain their labels. The edge $ox$ is isolated in the new support.

Compare the sets covered by the same sequence of labels before and after
the change. If an external base does not occur in a sequence that
originally covered all vertices, then both its spokes must occur.
After the change, all external vertices are still covered, and these
two spokes also supply $y$. To cover $x$, the sequence had to contain
$H$ or $X$, either of which now supplies $ox$. Thus every such sequence
still covers all vertices.

If all external bases occur, combine the external spoke labels into a
single label $S$ of weight $w_S$. The remaining conditions for complete
collection are
\[
\begin{array}{c|cc}
 & \text{before the change} & \text{after the change}\\ \hline
 Y\text{ occurs}                  & H\text{ or }X & H\text{ or }X\\
 X,Y\text{ both absent}           & H\text{ and }S & H\text{ and }S\\
 X\text{ occurs},\ Y\text{ absent}& H & S .
\end{array}
\]
In the last row, a sequence covers all vertices only before the change
if it contains $H$ and does not contain $S$. Replace all $j$ occurrences
of $H$ by $S$. This gives precisely the sequences that cover all vertices
only after the change. The inverse replaces every $S$ by $H$.
This replacement multiplies the probability of the sequence by
$(w_S/h)^j$, which is at least one.

Even when $w_S=h$, there is a strict gain from a sequence in which one
external base does not occur. Take $X$, both spokes $ou,ow$ of one
external petal, and the base of each of the other $\nu-2$ petals.
Originally $y$ was missing. After the change, the three spoke labels
supply $ox,yu,yw$. This uses $\nu+1$ draws. Appending any number of
copies of $X$ preserves the strict gain.

\paragraph{All bases are heavy.}
Each base has weight greater than the total weight of the spokes of
the other petals. First equalize the spoke weights in each petal while
preserving their sum. Fix the positions occupied by its spokes and the
outcomes of the other draws. If the base occurs, the choice of spokes
does not affect coverage. If it does not occur, both spokes are needed.
With $d\ge1$ positions for these spokes and conditional weights $z,1-z$,
the probability of seeing both is $1-z^d-(1-z)^d$, which is maximized at
$z=1/2$. For $d=0$, the weights do not matter. Thus equalizing the weights
does not decrease any $\Phi_k$.

Let $h_r$ be the base weight of petal $r$ and $v_r$ the weight of each
of its spokes. Choose two petals with the largest spoke weights $v_r$ and
label them $i,j$ so that $p=h_i+v_i\ge q=h_j+v_j$. In a single draw,
each private vertex of petal $i$ is covered with probability $p$ and
each private vertex of petal $j$ with probability $q$.
Set $x=v_i$, $y=v_j$, and $\tau=x+y$. Transfer weight $y$ from each
spoke of petal $j$ to the corresponding spoke of petal $i$.
Each private vertex of petal $i$ lies on one spoke and on the base.
To keep its probability of appearing equal to $p$, reduce the weight of
base $i$ by $y$. Increase the weight of base $j$ by $y$ to preserve the
probability $q$. After the transfer, only the base of petal $j$ remains
and forms an isolated edge.

With $p,q,\tau$ fixed, let $\eta$ denote the current weight of each
spoke of petal $i$. Then each spoke of petal $j$ has weight $\tau-\eta$,
and the base weights are determined by $p$ and $q$ as follows.
\[
\begin{array}{c|cc}
 \text{Petal} & \text{Weight of each spoke} & \text{Base weight}\\ \hline
 i & \eta      & p-\eta\\
 j & \tau-\eta & q-(\tau-\eta)
\end{array}
\]
Before the transfer $\eta=x$, and after it $\eta=\tau$. For
$x\le\eta\le\tau$, the mass transferred to each spoke of petal $i$ is
$\eta-x$. The two petals have total
weight $p+q+\tau$ for every $\eta$. The entire interval $0\le \eta\le\tau$
is feasible since $p\ge q$ and $q-\tau=h_j-x>0$ by the assumption that
base $j$ is heavy. Let $F_k(\eta)$ be the probability of complete collection within $k$ draws. We will show that $F_k(\tau)>F_k(x)$ for every $k\ge\nu+1$.

A word of length $n$ is a sequence of $n$ labels, and its weight is the
product of their weights. First consider all words of length $k$ that
cover all $2\nu$ private vertices of the bouquet. For complete collection,
such a word must also cover the center $o$, so it must contain at least
one spoke. To obtain $F_k(\eta)$, we therefore subtract from the total
weight of these words the weight of the words consisting entirely of
bases and containing every base at least once.

Call all petals other than $i$ and $j$ external. We use
inclusion--exclusion to account for coverage of the four private vertices
of petals $i$ and $j$. In each term, the allowed edges of these two petals
can be combined into a single label since none of them covers a private
vertex of any external petal. Call this label the free label and denote
the sum of the weights of the combined edges by $z$.

Let $G_n(z)$ be the total weight of the words of length $n$ that cover
all private vertices of the external petals when each position may be
occupied by any external edge or by the free label of weight $z$.
When all edges of $i$ and $j$ are allowed, $z=p+q+\tau$. To count words
without spokes, define $H_n(z)$ in the same way, allowing only bases
among the external edges. In this count, only the allowed bases of $i$
and $j$ contribute to $z$.

We now compute the weight of the words of length $k$ that cover all
private vertices. The function $G_k$ already accounts for coverage of
the external private vertices. Apply inclusion--exclusion to the four
private vertices of $i,j$. To forbid a vertex is to exclude every edge
containing it. Fix the external draws and their positions. Suppose that
$d$ positions remain for draws from $i,j$ and that the edges allowed
in these positions have weights $w_1,\ldots,w_m$. The total weight of
all ways of filling the positions is $(w_1+\cdots+w_m)^d=z^d$.
Thus, for any choice of forbidden vertices, it suffices to evaluate
$G_k$ at the total weight of the allowed edges of $i,j$.

Set $t=\tau-\eta$. If no vertex is forbidden in petal $i$, its allowed
edges have total weight $(p-\eta)+2\eta=p+\eta$. If one private vertex is
forbidden, only the spoke to the other private vertex remains, with
weight $\eta$. If both vertices are forbidden, no edge remains. The
corresponding weights in petal $j$ are $q+t$, $t$, and $0$.

The rows and columns of the table give the numbers of forbidden private
vertices in petals $i$ and $j$, respectively. A choice of $r$ forbidden
vertices in $i$ and $c$ in $j$ can be made in
$\binom{2}{r}\binom{2}{c}$ ways, and the sign of the term is $(-1)^{r+c}$.
The nine entries sum to the weight of the words covering all private
vertices.
\[
\begin{array}{c|ccc}
 & \multicolumn{3}{c}{\text{Forbidden vertices in petal }j}\\
 \text{Forbidden vertices in petal }i & 0 & 1 & 2\\ \hline
 0 & G_k(p+q+\tau) & -2G_k(p+\tau) & G_k(p+\eta)\\
 1 & -2G_k(q+\tau) & 4G_k(\tau)    & -2G_k(\eta)\\
 2 & G_k(q+t)      & -2G_k(t)      & G_k(0)
\end{array}
\]
For instance, when $r=1$ and $c=2$, only one spoke of petal $i$ is
allowed, and it has weight $\eta$. Three vertices are forbidden, and there
are two choices of the forbidden vertex in $i$. This gives the entry
$-2G_k(\eta)$. Since $\eta+t=\tau$, the parameter $\eta$ appears only in the
arguments $p+\eta$, $\eta$, $q+t$, and $t$.

From this sum subtract the weight of the words that do not cover the
center. Such words consist entirely of bases, and every base must
occur at least once to cover all private vertices. The external bases
are already accounted for in $H_k$. By inclusion--exclusion for the
two remaining bases, the weight to be subtracted is
\[
 H_k(p+q-\tau)-H_k(p-\eta)-H_k(q-t)+H_k(0).
\]
The first term allows both bases, whose total weight is
$(p-\eta)+(q-t)=p+q-\tau$. The second excludes base $j$, and the third
excludes base $i$. The last term restores the words containing neither
base, which were subtracted twice. Since the entire expression is
subtracted from the sum of the entries in the table, $H_k(p-\eta)$ and
$H_k(q-t)$ appear with positive signs in $F_k(\eta)$.

Let $C$ be the sum of the terms independent of $\eta$. The remaining terms for petal $i$ give
$G_k(p+\eta)-2G_k(\eta)+H_k(p-\eta)$, and those for petal $j$ give the same
expression with $q,t$ in place of $p,\eta$. Thus
\begin{equation}\label{s2:eq:bouquet-egf}
 \begin{aligned}
 F_k(\eta)&=C+g_p(\eta)+g_q(\tau-\eta),\\
 g_\ell(\eta)&=G_k(\ell+\eta)-2G_k(\eta)+H_k(\ell-\eta).
 \end{aligned}
\end{equation}

Fix $k\ge\nu+1$. Since $\tau=x+y$, equation~\eqref{s2:eq:bouquet-egf} gives
\[
 F_k(\tau)-F_k(x)
 =\bigl[g_p(x+y)-g_p(x)\bigr]-\bigl[g_q(y)-g_q(0)\bigr].
\]
Both increments are taken over intervals of length $y$. We will use a
bound for the external petals to prove that $g_p$ is strictly convex
and compare its increments over $[x,x+y]$ and $[0,y]$. We will then
compare the increments of $g_p$ and $g_q$ over the same interval $[0,y]$.

\emph{A bound for the external petals.}
Set $n=k-2$. Differentiating $J_k(z)$ twice with respect to $z$ marks
an ordered pair of distinct positions occupied by the free label.
Deleting them leaves a word of length $n$ subject to the same coverage
conditions. Hence
\[
 \frac{\mathrm d^2}{\mathrm dz^2}J_k(z)=k(k-1)J_n(z),
 \qquad J\in\{G,H\}.
\]
The polynomials $G_n,H_n$ have nonnegative coefficients and are
nondecreasing on $[0,\infty)$. By convexity,
$2G_n(\eta)\le G_n(0)+G_n(2\eta)$. Since $0\le\eta\le\tau\le p$,
equation~\eqref{s2:eq:bouquet-egf} gives
\[
 \frac{1}{k(k-1)}\frac{\mathrm d^2}{\mathrm d\eta^2}g_p(\eta)
 =G_n(p+\eta)-2G_n(\eta)+H_n(p-\eta)
 \ge H_n(p-\tau)-G_n(0).
\]
Thus, to prove that $g_p$ is strictly convex, it suffices to show that
$G_n(0)<H_n(p-\tau)$.

Consider a word counted by $G_n(0)$. It consists of external edges
and covers all private vertices of the external petals. For each external petal $r$ whose base
does not occur in the word, both of its spokes must occur. Replace
the first occurrence of its left spoke by its base and the first
occurrence of its right spoke by a special label indexed by $r$ and of
weight $v_r^2/h_r$. The product of the two weights remains $v_r^2$.
Replace each remaining spoke by an ordinary label of the same weight
that records its petal and side. Keep all original bases.

Every external base occurs in the resulting word. If the special label
indexed by $r$ occurs, base $r$ occurs exactly once. Replace this base
by the left spoke and the special label by the right spoke. The ordinary
labels identify the remaining spokes, and all other bases occupy their
original positions. This recovers the original word, so the map is
injective and preserves weight.

Set $V=\sum_{r\ne i,j}v_r$. The total weight of the ordinary and special
labels is
\[
 W=2V+\sum_{r\ne i,j}\frac{v_r^2}{h_r}.
\]
Fix all occurrences of the bases. Summing over the labels in each
remaining position gives a factor of $W$. Hence $H_n(W)$ is the total
weight of the words of length $n$ formed from these labels and bases
in which every external base occurs. This family contains the image
of the injection, so $G_n(0)\le H_n(W)$.

By the choice of $i,j$, we have $v_r\le\min(x,y)$ for every external
petal $r$. Since its base is heavy,
$h_r>2(x+y+V-v_r)\ge 2(y+V)$. Hence
\[
 \sum_{r\ne i,j}\frac{v_r^2}{h_r}
 \le \frac{y}{2(y+V)}\sum_{r\ne i,j}v_r
 =\frac{yV}{2(y+V)}<\frac y2.
\]
Since base $i$ is heavy, $p-\tau=h_i-y>2V+y$, and therefore $W<p-\tau$.
Construct a word of length $n$ by placing each of the $\nu-2$ external
bases once and filling the remaining $n-(\nu-2)=k-\nu\ge1$ positions
with the free label. Increasing the weight of this label from $W$ to
$p-\tau$ strictly increases the weight of this word, while the weights
of all other words do not decrease. When $\nu=2$, the word consists
entirely of free labels. Thus $G_n(0)\le H_n(W)<H_n(p-\tau)$.

\emph{Comparison before and after the transfer.}
The bound just proved gives strict convexity of $g_p$ on $[0,\tau]$. Since $x,y>0$, its
increment over $[x,x+y]$ exceeds its increment over $[0,y]$. We now
compare the increments of $g_p$ and $g_q$ over $[0,y]$. For
$q\le\ell\le p$ and $0\le \eta\le y$,
\[
 \frac{\partial^2}{\partial\ell\,\partial\eta}g_\ell(\eta)
 =k(k-1)\bigl[G_n(\ell+\eta)-H_n(\ell-\eta)\bigr]\ge0.
\]
Every word counted by $H_n(\ell-\eta)$ is also counted by $G_n(\ell+\eta)$
with at least the same weight, since $\ell+\eta\ge\ell-\eta$.
Thus the derivative of $g_\ell$ with respect to $\eta$ is nondecreasing
in $\ell$. Integrating with respect to $\eta$ from $0$ to $y$ gives the
second inequality below:
\[
 g_p(x+y)-g_p(x)>g_p(y)-g_p(0)\ge g_q(y)-g_q(0).
\]
The same transfer therefore gives $F_k(\tau)>F_k(x)$ for every
$k\ge\nu+1$, including when $p=q$. For $k<\nu+1$, it is impossible to
cover all $2\nu+1$ vertices.
\end{proof}

\subsection{Structure of optimal distributions}
\label{s2:sec:structure}

For odd $a\ge5$, we will prove that every optimal support consists of
isolated edges and one triangle. This description applies both to
minimizing the expected collection time and to maximizing the probability
of complete collection within $k$ draws for any fixed $k\ge\lceil a/2\rceil$.

To evaluate the same distribution $y$ at several horizons, we consider the
weighted sum $\sum_{j=0}^K c_j\Phi_j(y)$ with prescribed weights $c_j\ge0$.
In \S~7.7 we use $\lambda\Phi_i(y)+\mu\Phi_j(y)$ for two horizons
$i,j\ge\lceil a/2\rceil$ and positive weights $\lambda,\mu$.
By Theorem~\ref{s2:thm:structure}, a canonical distribution attains
the maximum of this sum. For odd $a\ge5$, this will allow us to choose
a single triangle mass such that, at both horizons, the probability of
complete collection under the resulting canonical distribution is at
least as high as under the original distribution.

The next lemma rules out leaves that do not belong to isolated edges.
We redistribute the weight of one original edge $uw$ among the pairs
containing $w$. Each draw still consists of two coupons.

\begin{lemma}[Ruling out leaves outside isolated edges]\label{s2:lem:singleton}
Let $a\ge4$, and suppose that in the support graph of $y$ the vertex
$v$ has a unique neighbor $u$, which has another neighbor $w$. Whenever
the pair $\{u,w\}$ is drawn, replace $u$ by a uniformly chosen coupon
from $[a]\setminus\{w\}$, independently on each such draw. Leave all
other pairs unchanged. The resulting distribution on pairs does not
decrease any $\Phi_k$ and strictly increases every $\Phi_k$ for which
$\Phi_k(y)>0$.
\end{lemma}

\begin{proof}
Every covering sequence contains $\{u,v\}$, since otherwise $v$ is
missing. This pair also retains $u$, the only coupon that can be lost
through the random replacement. Thus every previously covering sequence
still covers all coupons.

Fix $k$ with $\Phi_k(y)>0$ and take any covering sequence of $k$ draws.
Using pairs from the original support, we construct a sequence of the
same length that misses a coupon but can cover all coupons after the
random replacement. Choose an occurrence of an
edge $\{w,z\}$ containing $w$ in the original sequence and replace it
by $\{u,w\}$. If the sequence no longer covers all coupons, the random
replacement need only restore $\{w,z\}$ at this position.

If all coupons are still covered, also replace every occurrence of
$\{u,v\}$ by $\{u,w\}$. Now only $v$ is missing. At the first position
chosen, the random replacement can produce $\{w,v\}$, while at the
other modified positions it can retain $\{u,w\}$. At least one of
these other positions retains $u$, since $v$ was covered before the
replacements.

In either case, every other occurrence of $\{u,w\}$ can remain unchanged.
All the specified outcomes have positive probability and are independent.
Thus every sequence that originally covered all coupons still does so,
and with positive probability a sequence that originally missed a coupon
now covers them all. Hence $\Phi_k$ strictly increases.
\end{proof}

\begin{samepage}
\begin{theorem}[Optimal supports and weights]\label{s2:thm:structure}
Let $a\ge3$. Up to a permutation of the coupons, every maximizer of
$\Phi_k$ for $k\ge\lceil a/2\rceil$ and every minimizer of $\mathbb E T_y$
has the following form:
\begin{enumerate}
\item for $a=2q$, the uniform distribution on $q$ isolated edges;
\item for $a=3$, the uniform distribution on the three edges of a triangle;
\item for $a=2m+3\ge5$, $m$ isolated edges of equal weight and a triangle
whose three edges have equal weights.
\end{enumerate}
The same description holds for every maximizer of the weighted sum
$\sum_{j=0}^{K}c_j\Phi_j(y)$ if all $c_j\ge0$ and $c_j>0$ for at least
one $j\ge\lceil a/2\rceil$.
\end{theorem}
\end{samepage}

\begin{proof}
Set $h=\lceil a/2\rceil$. For $j<h$, at most $2j<a$ coupons can be
collected in $j$ draws, so $\Phi_j(y)=0$ for every $y$. If the weights
were positive only at such $j$, the weighted sum would be identically zero.

For a finite sum, set $c_j=0$ for $j>K$. Then
\[
 L_c(y):=\sum_{j\ge0}c_j\bigl(1-\Phi_j(y)\bigr)
 =\sum_{j=0}^K c_j-\sum_{j=0}^K c_j\Phi_j(y).
\]
The first term on the right is constant, so minimizing $L_c$ is
equivalent to maximizing the weighted sum. In particular, taking
$c_k=1$ and all other weights zero gives the problem of maximizing
$\Phi_k$. For the expected collection time, take $c_j=1$ for all
$j\ge0$. Then $L_c(y)=\mathbb E T_y$ by the tail-sum
formula~\eqref{pre:tail-equation}.

For a finite sum, a minimum exists by continuity on the simplex. For the
expected collection time, existence was established in
Proposition~\ref{pre:tail}. A distribution assigning positive weights to
all edges of a canonical support has \(\Phi_j>0\) for every \(j\ge h\)
and finite expected collection time. Thus for every minimizer there is
at least one index \(j_0\ge h\) such that \(c_{j_0}>0\) and
\(\Phi_{j_0}(y)>0\). Otherwise \(L_c=\sum_j c_j\) in the finite case and
the expected collection time is infinite in the second case. In
particular, the support graph of a minimizer has no isolated vertices.

A transformation that does not decrease any $\Phi_j$ does not increase
$L_c$. If it also strictly increases $\Phi_{j_0}$, then it strictly
decreases $L_c$.

For even \(a\), the assertion follows from Theorem~\ref{sdiv:thm:main}.
The uniform distribution on a perfect matching maximizes all \(\Phi_j\),
and for every \(j\ge h\) the inequality is strict for every distribution
of a different form. For \(a=3\) and \(j\ge1\), we have
\(\Phi_j=1-\sum_{i=1}^3 y_i^j\). Strict convexity of the power function
for \(j\ge2\) gives the same conclusion for the uniform distribution
on the triangle.

\textbf{The support for odd \(a\ge5\).}
Suppose that some minimizer of \(L_c\) has noncanonical support.
Among such minimizers, choose one with the largest possible number \(r\)
of isolated edges and, subject to this choice, the fewest positive-weight
edges. Set aside these \(r\) isolated edges. The remaining graph \(C\)
has an odd number \(n=a-2r\ge3\) of vertices, no isolated vertices, and
no \(K_2\) components. By our choice, every minimizer with more isolated
edges, or with the same number of isolated edges and fewer positive-weight
edges, has canonical support. Case~(iv) below produces one more isolated
edge, while case~(i) retains the original isolated edges and reduces
the number of positive-weight edges.

\begin{enumerate}
\item \textbf{A four-cycle.} If $C$ contains a four-cycle, apply
Lemma~\ref{s2:lem:cycle}. Passing to the tail probabilities $1-\Phi_j$
and summing the inequalities with weights $c_j$ yields
\[
 L_c(y)\ge(1-\lambda)L_c(y_-)+\lambda L_c(y_+).
\]
For the expected collection time, take the limit in the inequalities
for the partial sums, using the nonnegativity of all tail probabilities.
Since $y$ minimizes $L_c$ and $0<\lambda<1$, both endpoint
distributions $y_-,y_+$ are also optimal. Each endpoint retains the
original isolated edges and has fewer positive-weight edges. Thus both
supports are canonical. The strict part of the lemma now gives a strict
tail inequality for $j=j_0$. Since $c_{j_0}>0$, the inequality for $L_c$
is also strict, contradicting the optimality of the endpoints.
\item \textbf{A leaf.} If $v$ is a leaf with neighbor $u$, then $u$
has another neighbor since $C$ has no isolated edges. By
Lemma~\ref{s2:lem:singleton}, we can strictly increase $\Phi_{j_0}$
without decreasing the other probabilities. This contradicts the
minimality of $L_c$.
\end{enumerate}

Thus $C$ has no four-cycle and has minimum degree at least two. For
$n=3$ it is a triangle, and the entire support is canonical. Suppose
$n\ge5$. Any two distinct common neighbors of $u,v$ would form a
four-cycle with $u,v$. Thus if every two distinct vertices have a common
neighbor, it is unique. By Theorem~6 of Erd\H{o}s, R\'enyi, and
S\'os~\cite[p.~234]{ErdosRenyiSos1966}, in this case $C$ is a bouquet
of triangles. Two possibilities remain.
\begin{enumerate}[resume]
\item \textbf{A bouquet of triangles.} Keep the probability of drawing an edge of $C$ fixed and apply
Lemma~\ref{s2:lem:bouquet} to the conditional distribution given that
the drawn edge belongs to $C$.
Conditional on the positions of the draws from $C$ and the outcomes of
the other draws, the probability of complete collection does not
decrease. In $j_0$ draws, there is positive probability that each
external isolated edge occurs exactly once. This leaves
$j_0-r\ge(n+1)/2$ draws from $C$, for which the inequality is strict.
\item \textbf{Two vertices with disjoint neighborhoods.}
Choose two such vertices $P,Q$. Each has a neighbor outside $\{P,Q\}$
since the minimum degree is at least two. By Lemma~\ref{s2:lem:stars},
we can produce another isolated edge $PQ$ without decreasing any $\Phi_j$.
In particular, if $C$ has two triangle components, choose $P$ in one and
$Q$ in the other. On the four remaining vertices of the triangles the two
original edges survive, and each of the two neighbors of $P$ is joined
to each of the two neighbors of $Q$. Thus $K_3\sqcup K_3$ is transformed
into $K_2\sqcup K_4$.

The new distribution is optimal and has more isolated edges. By the choice
of the original minimizer, the new support must be canonical. This already
rules out two triangle components, since a canonical support has no
$K_4$ component.

In general, to apply the strict part of Lemma~\ref{s2:lem:stars}, we need
to cover $[a]\setminus\{P,Q\}$ in $h-1$ draws from the original
distribution. After deleting $P,Q$, every remaining vertex still has
an original neighbor. Indeed, vertices in $C$ lose at most one neighbor
each, and the original isolated edges survive. All remaining original
edges belong to the new canonical support. Thus each remaining isolated
edge is an original edge, and the original edges in the new triangle
cover all its vertices. The required cover consists of $m-1$ such
isolated edges and two original edges of the triangle, for a total of
$m+1=h-1$ draws, where $m=(a-3)/2$. Hence $\Phi_{j_0}$ strictly
increases, contradicting the minimality of $L_c$.
\end{enumerate}

It follows that the support of every minimizer consists of a triangle
and $m=(a-3)/2$ isolated edges.

\textbf{Weights.} Condition on which positions are occupied by isolated
edges and which by triangle edges. Equalizing the weights within each
group leaves the probability of this allocation unchanged. Complete
collection requires success in both groups, and the two events are
conditionally independent. For the isolated edges, equal weights are
optimal by Lemma~\ref{sdiv:lem:balance}. With $d\ge2$ draws from the
triangle, coverage fails only if the same edge is drawn all $d$ times.
For conditional edge weights $p_1,p_2,p_3$, this has probability
$\sum_i p_i^d$, uniquely minimized by equal weights. For $d<2$, the
triangle cannot be covered.

For strict inequality at $j_0$, it suffices to consider the event of
positive probability that exactly $m$ draws come from the isolated
edges and $d=j_0-m\ge2$ from the triangle. Given this event, the
probability of covering the isolated edges with conditional weights
$q_1,\ldots,q_m$ is $m!\prod_{i=1}^m q_i$, which is maximal only when
the weights are equal. Both conditional probabilities of success are
positive. Thus if the weights are unequal in either group, equalizing
them strictly increases $\Phi_{j_0}$, contrary to optimality.
\end{proof}

\subsection{Collection probabilities on the canonical support}
\label{s2:sec:profiles}

Let $m\ge1$. Each of the $m$ isolated edges has weight $u=(1-\theta)/m$,
and each edge of the triangle has weight $v=\theta/3$.
To collect all the coupons, we must draw every isolated edge and at least
two distinct edges of the triangle. Exclude $j$ specified isolated edges
and set $M=1-ju$. The total weight of the admissible sequences of length
$k$ that cover the triangle is
\begin{equation}\label{s2:eq:egf}
 M^k-3(M-2v)^k+2(M-3v)^k.
\end{equation}
Indeed, for each of the three triangle edges, we subtract the total weight
of the words containing neither of the other two triangle edges. Words containing no
triangle edge have thus been subtracted three times, so we add back twice
their weight. Inclusion--exclusion over the $m$ isolated edges
gives~\eqref{s2:eq:canonical}.

\begin{lemma}[Projection and concavity]\label{s2:lem:concavity}
Let $I_m^+=[L,R]$, where $L=3/(2m+3)$ and $R=2/(m+2)$.
For $k\ge m+2$,
\[
 \frac{\mathrm d}{\mathrm d\theta}F_{m,k}^{+}(\theta)>0
 \quad(0<\theta\le L),\qquad
 \frac{\mathrm d}{\mathrm d\theta}F_{m,k}^{+}(\theta)<0
 \quad(R<\theta<1).
\]
At $R$, the derivative is zero for $k=m+2$ and negative for $k\ge m+3$.
For $k=m+2$, the first derivative is positive throughout $0<\theta<R$.
The second derivative satisfies
\[
 \frac{\mathrm d^2}{\mathrm d\theta^2}F_{m,k}^{+}(\theta)<0
 \qquad (\theta\in I_m^+,\ k\ge m+2).
\]

For $k<m+2$, the probability $F_{m,k}^{+}$ is identically zero.
Replacing $\theta\in[0,1]$ by the nearest point of $I_m^+$ never
decreases $F_{m,k}^{+}$, for any $k$. If $\theta\notin I_m^+$
and $k\ge m+2$, this replacement strictly increases $F_{m,k}^{+}$.
\end{lemma}

\begin{proof}
If $k<m+2$, there are too few draws. Henceforth let $k\ge m+2$,
$0<\theta<1$, and $F(\theta)=F_{m,k}^+(\theta)$. Fix $r$ isolated edges and
$d$ edges of the triangle. Let $C_n(r,d)$ denote the total weight of the
words of length $n$ in which precisely these labels occur, each at least
once. We count the choices of the labels separately. The weight of a word
is the product of the weights of its labels, where $u=(1-\theta)/m$ and
$v=\theta/3$. For the empty word, $C_0(0,0)=1$.

\emph{The first derivative.}
Complete collection requires all isolated edges and either two or three
edges of the triangle. There are three ways to choose two of the three
edges, so $F=3C_k(m,2)+C_k(m,3)$.

Let $\mathcal E$ be the edge set of the canonical support and let
$\mathcal W_k\subseteq\mathcal E^k$ be the set of sequences
$\mathbf e=(e_1,\ldots,e_k)$ that cover all the coupons. Write
$p_e(\theta)$ for the weight of edge $e$. By independence, the weight of
a sequence is $\prod_{j=1}^k p_{e_j}$. Applying the product rule to the
weight of each covering sequence and collecting terms according to the
position $i$ of the differentiated factor gives
\[
\begin{aligned}
 \frac{\mathrm d}{\mathrm d\theta}F(\theta)
 &=\sum_{\mathbf e\in\mathcal W_k}
   \frac{\mathrm d}{\mathrm d\theta}\left(\prod_{j=1}^k p_{e_j}\right)
 =\sum_{\mathbf e\in\mathcal W_k}\sum_{i=1}^k
   \frac{\mathrm d p_{e_i}}{\mathrm d\theta}\prod_{j\ne i}p_{e_j}\\
 &=\sum_{i=1}^k
   \underbrace{\sum_{\mathbf e\in\mathcal W_k}
   \frac{\mathrm d p_{e_i}}{\mathrm d\theta}\prod_{j\ne i}p_{e_j}}
   _{\text{contribution of position }i}.
\end{aligned}
\]
For each fixed $i$, the sum in the last line runs over all covering
sequences. Interchanging positions $i$ and $k$ preserves coverage and
matches each term for position $i$ with a term of the same value for
position $k$. Thus all $k$ positions make the same contribution.

To compute the contribution of the last position, set $n=k-1$ and group
the sequences according to the first $n$ draws $z=(z_1,\ldots,z_n)$.
For fixed $z$, sum the derivatives of the weights of all edges $e$ for
which $(z,e)$ covers all the coupons. Then
\[
 \frac1k\frac{\mathrm d}{\mathrm d\theta}F(\theta)
 =\sum_{z\in\mathcal E^n}
  \left(\prod_{j=1}^n p_{z_j}\right)
  \sum_{\substack{e\in\mathcal E\\(z,e)\in\mathcal W_k}}
  \frac{\mathrm d p_e}{\mathrm d\theta}.
\]
For an isolated edge, $\frac{\mathrm du}{\mathrm d\theta}=-1/m$, and for
a triangle edge, $\frac{\mathrm dv}{\mathrm d\theta}=1/3$. If $z$ already
covers all the coupons, every edge is admissible and the inner sum is
$\frac{\mathrm d}{\mathrm d\theta}(mu+3v)=0$. If one further draw cannot
complete the collection, the inner sum is empty. This leaves two cases
in which the last draw completes the collection.

Suppose that the first $n$ draws contain all the isolated edges and just
one of the three triangle edges. There are three choices for this edge,
and for each choice the words have total weight $C_n(m,1)$. Either of
the other two edges completes the collection. The derivatives of their
weights sum to $2/3$, giving a contribution of
$3\cdot(2/3)C_n(m,1)=2C_n(m,1)$.

If the triangle is already covered and exactly one isolated edge is
missing, there are $m$ choices for the missing edge. For each choice the
words in the first $n$ positions have total weight
$3C_n(m-1,2)+C_n(m-1,3)$, since there are three ways to choose the two triangle edges that occur. Only the missing isolated edge completes the
collection. The derivative of its weight is $-1/m$, so summing over the
choices gives a factor of $m(-1/m)=-1$. Adding the two contributions yields
\begin{equation}\label{s2:eq:profile-first}
 \frac1k\frac{\mathrm d}{\mathrm d\theta}F(\theta)
 =2C_n(m,1)-\bigl[3C_n(m-1,2)+C_n(m-1,3)\bigr].
\end{equation}
In each of the three terms, fix $m-1$ isolated edges and one triangle edge
among the required labels. Suppose the remaining labels occupy exactly
$j\ge1$ positions. These positions can be chosen in $\binom nj$ ways.
The other $n-j$ positions must contain all the fixed labels and only those
labels. This gives a common factor of $\binom nj C_{n-j}(m-1,1)$.
In the first term a single additional isolated edge occupies all $j$
positions, giving weight $u^j$. In the second term a single additional
triangle edge occupies them, giving weight $v^j$. In the third term both
additional triangle edges must occur, so from $(2v)^j$ we subtract the
two words consisting of a single repeated edge, each of weight $v^j$.
Combining the coefficients gives
$2u^j-3v^j-\bigl((2v)^j-2v^j\bigr)=2u^j-(2v)^j-v^j$. Hence
\[
 \frac1k\frac{\mathrm d}{\mathrm d\theta}F(\theta)
 =\sum_{j=1}^{n}\binom nj C_{n-j}(m-1,1)
 \bigl[2u^j-(2v)^j-v^j\bigr].
\]
For $0<\theta<1$, the factor $C_{n-j}(m-1,1)$ is positive precisely when
$n-j\ge m$, so only $1\le j\le k-m-1$ need be considered.
The left endpoint $L$ is determined by $u=2v$. For $\theta\le L$ we
have $u\ge2v$, giving $2u^j-(2v)^j-v^j>0$ throughout this range of $j$.
The right endpoint $R$ is determined by $u=3v/2$. For $\theta\ge R$ we
have $u\le3v/2$, and $3v/2$ is the mean of $v$ and $2v$. By convexity
of the power function, $2u^j\le2(3v/2)^j\le v^j+(2v)^j$, with strict
inequality in the second step for $j\ge2$. When $k=m+2$, the only
bracket that remains is $2u-3v$, which is positive to the left of $R$,
zero at $R$, and negative to its right. When $k\ge m+3$, the contribution
for $j=2$ is strictly negative at $R$. To the right of $R$, the bracket
for $j=1$ is also strictly negative. These signs and the continuity
of $F$ on $[0,1]$ prove the assertion about projection.

\emph{The second derivative.}
For $k\ge m+2$ we show that $F$ is strictly concave on $I_m^+$.
Differentiate~\eqref{s2:eq:profile-first} with respect to $\theta$.
The weight of each word on the right is a product of $k-1$ label
weights, so the product rule gives one term for each of the $k-1$
positions. After summation over the words, all positions make the same
contribution, since permuting the positions preserves the set of labels
that occur. Taking into account the factor $1/k$ on the left
of~\eqref{s2:eq:profile-first}, we see that the contribution of the last
position is $\frac{1}{k(k-1)}\frac{\mathrm d^2}{\mathrm d\theta^2}F(\theta)$.

Each row of the table specifies a condition on the word of length $k-1$.
The first row requires all $m$ isolated edges and exactly one specified
triangle edge. Its coefficient in~\eqref{s2:eq:profile-first} is $2$.
The second row requires exactly $m-1$ specified isolated edges and
either two or three triangle edges. Both possibilities cover the triangle. The minus sign in front of the
bracket in~\eqref{s2:eq:profile-first} applies to both.

Fix the first $k-2$ draws. The columns indicate what is still needed to
satisfy the condition in a given row. If the condition already holds,
the last label must preserve it. In the other two columns, one required
isolated edge or one triangle edge is still needed, respectively.
If no single label can satisfy the condition, the contribution is zero.

Write $C(r,d)=C_{k-2}(r,d)$. Then $P_r=3C(r,2)+C(r,3)$ is the total
weight of the words of length $k-2$ that cover the triangle and contain
precisely $r$ specified isolated edges. To obtain each table entry,
multiply the weight of each word in the corresponding column by the row
coefficient and by the sum of the derivatives of the weights of the
admissible final labels, then sum over these words.
\begin{center}
\begingroup
\small
\setlength{\tabcolsep}{4pt}
\begin{tabular}{@{}lccc@{}}
\toprule
Row condition & \shortstack{Already\\satisfied}
 & \shortstack{One isolated\\edge missing}
 & \shortstack{Triangle\\edge needed} \\
\midrule
\shortstack[l]{$m$ isolated edges and\\exactly one specified triangle edge}
 & $-\tfrac43 C(m,1)$ & $-2C(m-1,1)$ & $\tfrac23 C(m,0)$ \\[3pt]
\shortstack[l]{$m-1$ specified isolated edges\\and the triangle covered}
 & $-P_{m-1}/m$ & $\tfrac{m-1}{m}P_{m-2}$ & $-2C(m-1,1)$ \\
\bottomrule
\end{tabular}
\endgroup
\end{center}

Consider the bottom right entry. All $m-1$ specified isolated edges have
occurred in the first $k-2$ draws, together with exactly one triangle
edge. There are three choices for this edge, so the words have total
weight $3C(m-1,1)$. Either of the other two triangle edges is admissible,
and the derivatives of their weights sum to
$2\,\frac{\mathrm dv}{\mathrm d\theta}=\tfrac23$. Including the row
coefficient gives
\[
 (-1)\cdot 3C(m-1,1)\cdot\tfrac23=-2C(m-1,1).
\]

In the first column, the two triangle edges other than the specified one
are excluded in the first row, and one isolated edge is excluded in the
second. The total weights of the admissible labels are $1-2v$ and $1-u$,
respectively, with derivatives $-2/3$ and $1/m$ with respect to $\theta$.
In the middle column the missing isolated edge can be chosen in $m$ and
$m-1$ ways, respectively, and the derivative of its weight is $-1/m$.
When $m=1$, the bottom entry in this column is zero. In the top right
entry only the specified triangle edge is admissible, giving
$2\cdot\tfrac13 C(m,0)$. Summing the entries gives
\begin{equation}\label{s2:eq:profile-second}
 \begin{aligned}
 \frac{1}{k(k-1)}\frac{\mathrm d^2}{\mathrm d\theta^2}F(\theta)
 ={}&\frac23 C(m,0)-\frac43 C(m,1)-4C(m-1,1)\\
 &-\frac1m P_{m-1}+\frac{m-1}{m}P_{m-2}.
 \end{aligned}
\end{equation}
When $m=1$, the last term is absent.

We now show that the right-hand side of~\eqref{s2:eq:profile-second} is
negative on $I_m^+=[L,R]$. On this interval $3v/2\le u\le2v$.
Fix one of the $m$ isolated edges and two triangle edges. In a word
counted by $C(m,0)$, replace each of the $t\ge1$ occurrences of the fixed
isolated edge by either of the chosen triangle edges. Summing over all
replacements changes the factor $u^t$ to $(2v)^t\ge u^t$.

If only one of the two triangle edges occurs, the resulting words have
total weight $A:=C(m-1,1)$ for each of the two choices of that edge.
If both occur, the total weight is $B:=C(m-1,2)$.
In either case the other labels are the $m-1$ specified isolated edges.
Replacing every triangle edge by the fixed isolated edge recovers the
original word, so each resulting word is counted exactly once.
For an original word $\omega$, let $t_\omega$ be the number of
occurrences of the label being replaced and let $w_\omega$ be the
product of the weights in the remaining positions. Since $u\le 2v$,
summing over the original words gives
\[
 C(m,0)=\sum_\omega w_\omega u^{t_\omega}
 \le \sum_\omega w_\omega(2v)^{t_\omega}=2A+B.
\]
At $\theta=L$ we have $u=2v$, and equality holds.

Now consider the words counted by $B=C(m-1,2)$. Fix one of the two
specified triangle edges and replace all its occurrences by the missing
isolated edge. The image consists precisely of the words in which all $m$ isolated
edges occur and the only triangle edge that occurs is the other specified edge. The replacement is reversible. In each of the
$t_\omega$ positions where a replacement is made, the weight changes
from $v$ to $u\ge v$. Thus summing over the original words counted by
$B$ gives
\[
 C(m,1)=\sum_\omega w_\omega u^{t_\omega}
 \ge \sum_\omega w_\omega v^{t_\omega}=B.
\]

For $m\ge2$, apply the same replacement to the words counted by
$C(m-2,2)$, choosing one of the missing isolated edges in advance.
The resulting words contain $m-1$ specified isolated edges and one
specified triangle edge. The weight of each word is multiplied by
$(u/v)^{t_\omega}\ge 3/2$, since $u/v\ge 3/2$ and the edge being
replaced occurs at least once. Hence
\[
 C(m-2,2)=\sum_\omega w_\omega v^{t_\omega}
 \le \tfrac23\sum_\omega w_\omega u^{t_\omega}
 =\tfrac23 C(m-1,1)=\tfrac23 A.
\]
For the words counted by $C(m-2,3)$, two specified triangle edges remain
after this replacement. The same weight comparison gives
$C(m-2,3)\le\tfrac23 C(m-1,2)=\tfrac23 B$.
Multiplying the first bound by $3$ and adding the second yields
$P_{m-2}\le\tfrac23(3A+B)$ for $m\ge2$.
Since $(m-1)/m\le1$, the sum of the terms with positive coefficients
in~\eqref{s2:eq:profile-second} is at most
$\tfrac23(2A+B)+\tfrac23(3A+B)=\tfrac{10}3 A+\tfrac43 B$ for all $m\ge1$.
Dropping the nonpositive term $-P_{m-1}/m$ and using $C(m,1)\ge B$, we obtain
\[
\begin{aligned}
 \frac{1}{k(k-1)}\frac{\mathrm d^2}{\mathrm d\theta^2}F(\theta)
 &\le \tfrac{10}3 A+\tfrac43 B-4A-\tfrac43 C(m,1)\\
 &\le -\tfrac23 A<0,\qquad \theta\in I_m^+ .
\end{aligned}
\]
Indeed, $A=C_{k-2}(m-1,1)>0$. Each of the $m-1$ specified isolated edges
and the specified triangle edge has positive weight. Since $k-2\ge m$,
we may use each of these labels once and fill the remaining positions
with any one of them.
\end{proof}

\subsection{Maximum probability of complete collection}

\begin{proof}[Proof of Theorem~\ref{s2:thm:horizon}]
For $k\ge\lceil a/2\rceil$, Theorem~\ref{s2:thm:structure} characterizes
all optimal supports and shows that the triangle edges have equal weights,
as do the isolated edges. For $a=2q$, this leaves the ordinary coupon
collector problem with $q$ equally likely types, and inclusion--exclusion
gives~\eqref{s2:eq:evenvalue}. For $a=3$, the collection is complete once at least two distinct
triangle edges have been drawn. Under the uniform distribution on the
three edges, this occurs within $k$ draws with probability $1-3^{1-k}$
for $k\ge1$. For $k<\lceil a/2\rceil$, complete collection
is impossible under any distribution.

Suppose $a=2m+3\ge5$ and $k\ge m+2$. By Lemma~\ref{s2:lem:concavity},
the maximum is attained on $I_m^+=[L,R]$. On this interval the
derivative strictly decreases from a positive value at $L$ to a
nonpositive value at $R$, and it does not vanish on
$(0,1)\setminus[L,R]$. Thus the maximum is attained at the unique critical point
$\theta_{m,k}^{+}\in(0,1)$. This point equals $R$
when $k=m+2$ and lies in $(L,R)$ when $k\ge m+3$.

For $k=m+2$, all $m$ isolated edges and two distinct triangle edges must
occur. The two triangle edges can be chosen in three ways, and the
resulting $m+2$ draws can be ordered in $(m+2)!$ ways. Hence
\begin{equation}\label{s2:eq:firstprofiles}
 F_{m,m+2}^{+}(\theta)=(m+2)!\,u^m\,3v^2
 =\frac{(m+2)!}{3m^m}\theta^2(1-\theta)^m.
\end{equation}
At the maximum $\theta=2/(m+2)$, each isolated edge has mass $1/k$
and each triangle edge has mass $2/(3k)$. Substitution gives
\begin{equation}\label{s2:eq:firstvalue}
 \theta_{m,m+2}^{+}=\frac{2}{m+2},\qquad
 \alpha_{\max}(2k-1,2k-3,k)=\frac{4k!}{3k^k}\quad(k\ge2).
\end{equation}
When $k=2$, the second formula gives $2/3$, the probability of
collecting three coupons in two draws under the uniform distribution
on the three edges. For $k\ge m+3$ we have
$\theta_{m,k}^{+}<R=2/(m+2)$, so the optimal mass of the triangle
depends on the number of draws.
\end{proof}

\subsection{Minimum expected collection time}

\begin{proof}[Proof of Theorem~\ref{s2:thm:mean}]
By Theorem~\ref{s2:thm:structure}, all minimizers belong to the canonical
family. For $a=2q$, Theorem~\ref{sdiv:cor:mean} with $s=2$ gives
$M(2q,2q-2)=qH_q$. For $a=3$, after the first draw it remains to wait
for a different triangle edge. Each draw gives such an edge with
probability $2/3$, so $M(3,1)=1+3/2=5/2$.

Suppose $a=2m+3\ge5$. At $\theta=0,1$ complete collection is
impossible. Let $0<\theta<1$, $u=(1-\theta)/m$, and $v=\theta/3$.
Let $S$ be the number of draws needed to see all $m$ isolated edges.
With $j$ isolated edges still missing, the next draw yields one of them
with probability $ju$. The mean waiting time for the next new isolated edge
is therefore $1/(ju)$ draws. Summing over the numbers of missing edges gives
\[
 \mathbb E S=\sum_{j=1}^{m}\frac{1}{ju}=\frac{H_m}{u}.
\]

Let $T$ be the collection time. By time $S$ all isolated edges have
occurred. If at least two distinct triangle edges have occurred, then
$T=S$. If only one of the triangle edges has occurred, each subsequent
draw gives one of the other two with probability $2v$, so the expected
waiting time is $1/(2v)$ draws. If no triangle edge has occurred, the
expected waiting time for the first one is $1/\theta$ draws, followed
by a further $1/(2v)$ draws on average for a different triangle edge.
Both $S$ and the first time a triangle edge is drawn are stopping times.
After either time, subsequent draws are independent of the draws up to that
time and have the original distribution.

Let $P_c$ be the probability that by time $S$ none of the edges in a
specified set of triangle edges of total mass $c>0$ has occurred.
For this event to occur, at each stage the next new isolated edge must
precede any edge in the specified set. With $j$ isolated edges still
missing, the probability of this is
$ju/(ju+c)$, since all other draws can be ignored. Considering
$j=m,\ldots,1$ in turn gives
\begin{equation}\label{s2:eq:mean-sum}
 P_c=\prod_{j=1}^m\frac{ju}{ju+c}.
\end{equation}

The probability that no triangle edge has occurred by time $S$ is
$P_\theta$. To find the probability that fewer than two distinct
triangle edges have occurred, sum over the three pairs of triangle
edges the probability that both edges of the pair are absent.
If precisely one edge has occurred, there is a unique absent pair.
If none has occurred, all three pairs are absent. Thus the required
probability is $3P_{2v}-2P_\theta$. Using the expected waiting time at
each stage gives
\[
 \mathbb E(T-S)=\frac{P_\theta}{\theta}
 +\frac{3P_{2v}-2P_\theta}{2v}
 =\frac{3P_{2v}}{2v}-\frac{2P_\theta}{\theta}.
\]
Adding $\mathbb E S$ and substituting~\eqref{s2:eq:mean-sum} gives
\begin{equation}\label{s2:eq:mean-product}
 E_m^+(\theta)=\frac{H_m}{u}
 +\frac{9}{2\theta}\prod_{j=1}^m\frac{ju}{ju+2\theta/3}
 -\frac2\theta\prod_{j=1}^m\frac{ju}{ju+\theta},
 \qquad u=\frac{1-\theta}{m}.
\end{equation}

Thus $E_m^+$ is finite and smooth on $(0,1)$. By the tail-sum
formula~\eqref{pre:tail-equation},
\[
 E_m^+(\theta)=\sum_{k\ge0}\bigl(1-F_{m,k}^+(\theta)\bigr).
\]
The summands are constant for $k<m+2$ and strictly convex on $I_m^+$
for $k\ge m+2$ by Lemma~\ref{s2:lem:concavity}. Therefore $E_m^+$ is
strictly convex on this interval.

Write $I_m^+=[L,R]$. By Lemma~\ref{s2:lem:concavity}, each tail
probability $1-F_{m,k}^+$ is nonincreasing on $(0,L]$ and nondecreasing
on $[R,1)$. For $k=m+3$ its derivative is negative on the first interval
and positive on the second, including $R$, where the derivative for
$k=m+2$ is zero. The sum of the other tail probabilities is smooth,
being the difference between $E_m^+$ and the term with $k=m+3$.
This sum is nonincreasing on $(0,L]$ and nondecreasing on $[R,1)$. Thus the derivative of
$E_m^+$ is negative on $(0,L]$ and positive on $[R,1)$. By strict
convexity of $E_m^+$ on $[L,R]$, its derivative vanishes at exactly one
point $\theta_{m,E}^+\in(L,R)$, the unique minimizer.
\end{proof}

\subsection{Simultaneous comparison}

\begin{proof}[Proof of Theorem~\ref{s2:thm:dom}]
For even $a$, the uniform distribution on a perfect matching maximizes
every $\Phi_k$. For $a=3$, the uniform distribution on the triangle edges
has the same property (Theorem~\ref{s2:thm:structure}). In both cases
this single distribution works for all $k$.

Suppose $a=2m+3\ge5$. Fix an arbitrary distribution $y$ on pairs. Denote by
$f_k(\theta)=F_{m,k}^+(\theta)$ the probability of collecting all the
coupons within $k$ draws under the canonical distribution that assigns total
mass $\theta$ to the triangle. Set
\[
 J_k=\{\theta\in I_m^+: f_k(\theta)\ge\Phi_k(y)\}.
\]
We seek a single mass belonging to every $J_k$. If $k<m+2$, complete
collection is impossible, so $J_k=I_m^+$. For $k\ge m+2$, Theorem~\ref{s2:thm:structure} gives a canonical
distribution maximizing $\Phi_k$. By Lemma~\ref{s2:lem:concavity},
we may choose its triangle mass in $I_m^+$. Thus $J_k$ is nonempty.
It is closed by continuity of $f_k$, and by concavity it contains the
entire interval between any two of its points. Write $J_k=[\ell_k,r_k]$,
where $\ell_k$ and $r_k$ denote the left and right endpoints of this
nonempty closed interval.

We show that any two intervals $J_i$ and $J_j$ intersect. Suppose
otherwise and choose the indices so that $J_i$ lies to the left of $J_j$.
There is then a gap between them, with $r_i<\ell_j$, since equality
would give a common endpoint. Both $i$ and $j$ are at least $m+2$,
since for smaller indices the corresponding interval is all of $I_m^+$.

For $k\in\{i,j\}$, put $g_k(\theta)=f_k(\theta)-\Phi_k(y)$.
For $\theta\in I_m^+$, this difference is nonnegative precisely when $\theta\in J_k$. Thus both
differences are negative in $(r_i,\ell_j)$, and at the endpoints
\[
 g_i(r_i)=0,\quad g_j(r_i)<0,\qquad
 g_i(\ell_j)<0,\quad g_j(\ell_j)=0.
\]
Indeed, if $g_i(r_i)>0$, then continuity would imply that $J_i$
contains points to the right of its right endpoint $r_i$. The same
argument shows that $g_j$ vanishes at its left endpoint $\ell_j$.

Choose positive weights so that the weighted difference
$g=\lambda g_i+\mu g_j$ takes equal values at $r_i$ and $\ell_j$.
At the former point only $\mu g_j(r_i)$ remains, and at the latter only
$\lambda g_i(\ell_j)$ remains. It therefore suffices to take
$\lambda=-g_j(r_i)$ and $\mu=-g_i(\ell_j)$. Then
$g(r_i)=g(\ell_j)=-\lambda\mu<0$. Throughout $(r_i,\ell_j)$ we have $g<0$,
since both terms are negative.

We claim that $g<0$ on the rest of $I_m^+$ as well. The function $g$ is
concave, being a positively weighted sum of concave functions.
If $\theta<r_i$, then $r_i$ lies between $\theta$ and $\ell_j$.
Writing $r_i=t\theta+(1-t)\ell_j$ with $0<t<1$, we obtain
\[
 g(r_i)\ge t g(\theta)+(1-t)g(\ell_j)
        =t g(\theta)+(1-t)g(r_i).
\]
It follows that $g(\theta)\le g(r_i)<0$. If $\theta>\ell_j$, then
$\ell_j$ lies between $r_i$ and $\theta$, and the same argument gives
$g(\theta)\le g(\ell_j)<0$. Thus for every $\theta\in I_m^+$,
\[
 \lambda f_i(\theta)+\mu f_j(\theta)
 <\lambda\Phi_i(y)+\mu\Phi_j(y).
\]
By Theorem~\ref{s2:thm:structure}, the maximum of
$\lambda\Phi_i+\mu\Phi_j$ over all distributions is attained by a
canonical distribution. By Lemma~\ref{s2:lem:concavity}, its parameter
may be chosen in $I_m^+$. At that parameter the left-hand side is the
maximum itself, and the right-hand side is the value of the same sum at $y$. This contradiction proves that
$J_i\cap J_j\ne\varnothing$.

It remains to choose one triangle mass that works for every $k$. Pairwise
intersection implies $\ell_i\le r_j$ for all $i,j$, since otherwise
the whole of $J_i$ would lie to the right of $J_j$. Thus every right
endpoint $r_j$ is an upper bound for all the left endpoints, so
\[
 \ell_j\le\theta^*\coloneqq\sup_{i\ge0}\ell_i\le r_j
 \qquad (j\ge0).
\]
Hence $\theta^*$ belongs to every $J_k$, and the canonical
distribution $y^*$ with triangle mass $\theta^*$ satisfies
$\Phi_k(y^*)\ge\Phi_k(y)$ simultaneously for all $k\ge0$.
\end{proof}

\section{Further questions}

The first cases in which our results do not fully characterize the optimal
distributions occur with seven coupons, for $(a,b)=(7,3)$ and $(7,4)$.
In the former case, three coupons are missing from each draw. In the
latter, three coupons are drawn, so $s=3$. In both cases it remains to determine the distributions that
maximize the probability of collecting all coupons within $k$ draws for
$k\ge3$ and those that minimize the expected collection time.

A partition of the seven coupons into a set of three and a set of four
gives a natural family of candidates. The set of three is chosen with
probability $1-4p$, and each of the four triples within the set of four is
chosen with probability $p$, where $0<p<1/4$. When $b=3$, the chosen
triple is the missing set, and when $s=3$, it is the drawn set. For
$(a,b)=(7,4)$, Schilling's Conjecture~6.2.2~\cite[p.~73]{Schilling2020}
predicts that a distribution in this family attains the minimum expected
collection time. Does a simultaneous comparison as in
Theorem~\ref{s2:thm:dom} hold in each of these two cases? In other words,
given any distribution $y$, can we choose a single distribution $y^*$
from this family such that $\Phi_k(y^*)\ge\Phi_k(y)$ for all $k\ge0$?

\end{document}